\documentclass{article}
\title{}
\date{}
\author{}
\usepackage{amssymb}
\usepackage{amsfonts}
\usepackage{mathrsfs}
\usepackage{amsmath}
\usepackage{graphicx}
\usepackage{empheq}
\usepackage{indentfirst}
\usepackage{cases}
\usepackage[titletoc]{appendix}
\usepackage{hyperref}
\newcommand{\cqfd}
{
\mbox{}%
\nolinebreak%
\hfill%
\rule{2mm}{2mm}%
\medbreak%
\par%
}

\newtheorem{thm}{\indent Theorem}[section]
 
 \newtheorem{lemma}{\indent Lemma}[section]
 \newtheorem{proposition}{\indent Proposition}[section]
 
 \newtheorem{rem}{\indent Remark}[section]
 \numberwithin{equation}{section}
\allowdisplaybreaks

\begin{document}
\title{Boundary exponential stabilization of 1-D inhomogeneous quasilinear hyperbolic systems}
\author{Long Hu\thanks{L. Hu is with School of Mathematics, Shandong University, Jinan, Shandong 250100, China. Sorbonne Universit\'{e}s, UPMC Univ Paris 06, UMR 7598, Laboratoire
Jacques-Louis Lions, F-75005, Paris, France. School of Mathematical
Sciences, Fudan University, Shanghai 200433,
China. E-mail: \texttt{hu@ann.jussieu.fr},\ \ \texttt{hul10@fudan.edu.cn}. This author
was supported by the China Scholarship Council for Ph.D. study at UPMC (No. 201306100081) and was partially supported by ERC advanced grant 266907
(CPDENL) of the 7th Research Framework
Programme (FP7). },  Rafael Vazquez\thanks{Rafael Vazquez is with the Department of Aerospace Engineering, Universidad de Sevilla, Camino de los Descubrimiento s.n., 41092 Sevilla, Spain.
        {\texttt{rvazquez1@us.es}}}, Florent Di Meglio\thanks{F. Di Meglio is with MINES ParisTech, PSL Research University, CAS - Centre automatique et syst\`emes, 60 bd St Michel, 75006 Paris, France.
        {\texttt{florent.di\_meglio{\@}mines-paristech.fr}}}\ \ and\ Miroslav Krstic\thanks{M. Krstic is with the Department of Mechanical and Aerospace Engineering, University of California San Diego, La Jolla, CA 92093-0411, USA.
        {\texttt{krstic@ucsd.edu}}}}

\maketitle
\begin{abstract}
This paper deals with the problem of boundary stabilization of first-order $n\times n$ inhomogeneous quasilinear  hyperbolic systems. A backstepping method is developed.  
The main result supplements the previous works on how to design multi-boundary feedback controllers to realize exponential stability of the original nonlinear system in  the spatial $H^2$ sense.  
\end{abstract}



\section{Introduction and Main Result}\label{problem}
Consider the following  1-D $n\times n$ inhomogeneous quasilinear hyperbolic system
\begin{align}
\frac{\partial u}{\partial t}+A(x,u) \frac{\partial u}{\partial x}=F(x,u),\ x\in[0,1],\ t\in[0,+\infty),\label{3.1}
\end{align}
where, $u=(u_{1},\ldots,u_{n})^{T}$is an unknown vector function of $(t,x)$,\
$A(x,u)$ is an $n\times n$ matrix with $C^2$ entries $a_{ij}(x,u)(i,j=1,\cdots,n)$, $F:[0,1]\times\mathbb{R}^n  \rightarrow\mathbb{R}^n$ is a vector valued function with $C^2$ components $f_i(x,u)(i=1,\cdots,n)$ with respect to $u$ and
\begin{align}
F(x,0)\equiv 0.
\end{align}
Denote 
\begin{align}
&\frac{\partial F}{\partial u}(x,0):=(f_{ij}(x))_{n\times n},
\end{align}
we assume that $f_{ij}\in C^2([0,1])$

By the definition of hyperbolicity,  we assume that $A(x,0)$ is a diagonal matrix with distinct and nonzero eigenvalues $A(x,0)=\mathrm{diag}(\Lambda_1(x),\cdots,\Lambda_n(x))$, which are, without loss of generality,  ordered as follows:
\begin{align}\label{3.7}
\Lambda_1(x) < \Lambda_2(x) < \cdots < \Lambda_m(x) < 0 
< \Lambda_{m+1}(x) <\cdots<  \Lambda_n(x),\forall x\in[0,1].
\end{align}
Here  and in what follows, $\mathrm{diag}(\Lambda_1(x),\cdots,\Lambda_n(x))$ denotes the diagonal matrix whose $\mathrm{i}$-th element on the diagonal is $\Lambda_i(x)$.


Under the assumption (\ref{3.7}),  a general kind of boundary conditions which guarantee the well- posedness of the forward problem on the domain $\{(t, x) | t \geq 0, 0\leq x \leq 1\}$ can be written as (see \cite{Li-1994}): 
\begin{align}
&x=0: u_s=G_s(u_1,\cdots,u_m), \ \ s=m+1,\cdots,n,\label{3.10}\\
&x=1: u_r=h_r(t), \ \ r=1,\cdots,m,\label{3.11}
\end{align}
where $G_s$  are $C^2$ functions, and  we assume that they vanish at the origin, i.e. 
\begin{align}\label{3.12}
G_s(0,\cdots,0)\equiv 0, s=m+1,\cdots,n,
\end{align}
while $H=(h_1,\cdots,h_m)^T$ are boundary controls. Our concern, in this paper, is to design a feedback control law for $H(t)$ in order to ensure that the closed-loop system is locally exponentially stable in the $H^2$ norm. 

In other words, we are interested in the following stabilization problem for the system (\ref{3.1}) and (\ref{3.10})-(\ref{3.11}):

\textbf{Problem (ES).} For any given $\lambda>0$. Suppose that $C^1$ compatibility conditions   are satisfied at the point $(t, x) = (0, 0)$. Does there exist a linear feedback control $\mathcal{B}: (H^2(0,1))^n\rightarrow \mathbb{R}^m$, verifying the $C^1$ compatibility conditions at the point $(t,x)=(0,1)$, such that for some $\varepsilon>0$,  the mixed initial-boundary  value problem (\ref{3.1}), (\ref{3.10})-(\ref{3.11})   and the initial conditions
\begin{align}\label{initial}
t=0: u(0,x):=\phi(x)=(\phi_1(x),\cdots,\phi_n(x)),
\end{align}
with $H(t)=\mathcal{B}(u(t,\cdot))$ admits a unique $C^0([0,\infty);(H^2(0,1))^n)$ solution $u=u(t,x)$, which satisfies
\begin{align}
\|u(t,\cdot)\|_{H^2(0,L)}\leq C e^{-\lambda t}\|\phi(\cdot)\|_{H^2(0,L)},
\end{align}
for some $C>0$, provided that $\|\phi(\cdot)\|_{H^2(0,L)}\leq \varepsilon$? 


The boundary stabilization problem for linear and nonlinear hyperbolic system has been widely studied in the last three decades or so. During this time, three parallel mathematical approaches have emerged. The first one is  the so-called ``Characteristic method", i.e. computing corresponding bounds by using explicit evolution of the solution along the characteristic curves. With this method, Problem (ES) has been previously investigated by Greenberg and Li (see \cite{Greenberg-Li-1984}) for $2\times 2$ systems and Li and Qin (see \cite{Li-1994,Qin-1985}) for a generalization to $n\times n$ homogeneous systems in the framework of $C^1$ norm.  Also, this method was developed by Li and Rao \cite{Li-Rao-2003} to study the exact boundary controllability for general inhomogeneous quasilinear hyperbolic systems.

The second method is the ``Control Lyapunov Functions method'', which is a useful tool to analyze the asymptotic behavior of dynamical systems. This method was first used by Coron et.al. to  design  dissipative boundary conditions for nonlinear homogeneous hyperbolic systems in the context of both $C^1$ and $H^2$ norm \cite{JMC-Bastin-2014,JMC-Bastin-Andera-2008,Coron-Andrea-Bastin-2007}. More recently, it has been shown in~\cite{JMC-Nguyen-2014} that the exponential stability strongly depends on the considered norm, i.e. a previously known sufficient condition for exponential stability with respect to the $H^2$ norm is not sufficient in the framework of $C^1$ norm. Although the Control Lyapunov Functions method has been introduced to study exponential stability for hyperbolic systems of balance laws, however, finding  a ``good'' Lyapunov Function is the main difficulty, especially  when the ``natural'' control Lyapunov functions do not lead to arbitrarily large exponential decay rate to the original system (see \cite{Bastin-Coron-2011-SCL}, \cite[Pages 314 and 361--371]{Coron-2007}).  This phenomenon indeed happens when we deal with Problem (ES) for the inhomogeneous hyperbolic systems (see \cite{JMC-Bastin-Andera-2008} and \cite{Coron-Andrea-Bastin-2007}). 

The third one is the ``Backstepping method'', which is now a popular mathematical tool to stabilize the finite dimensional and infinite dimensional dynamic systems (see \cite{Krstic-Smyshlyaev-2008, Krstic-2008-book, Smyshlyaev-Krstic-2010-book,Vazquez-Krstic-2008,Vazquez-Krstic-2-2008}).  In~\cite{JMC-Vazquez}, a full-state feedback control law, with actuation on only one end of the domain, which achieves $H^2$ exponential stability of the closed-loop $2\times 2$ linear and quasilinear hyperbolic system is derived using a backstepping method. Moreover, this method ensures that the linear hyperbolic system vanishes in finite time. Unfortunately, the method presented in \cite{JMC-Vazquez} can not be directly extended to $n\times n$ cases, especially when several states convecting in the same direction are controlled (see also \cite{Meglio}). In~\cite{hu-meglio}, a first step towards generalization to $3\times 3$ linear hyperbolic systems is addressed, in the case where  two controlled states are considered. 
With a similar Volterra transformation, designing an appropriate form of the target system, Hu et.al.\cite{hu-dimeglio-vazquez-krstic} adopt a classical backstepping controller to handle the Problem (ES) for general $n\times n$ linear hyperbolic systems. Well-posedness of the system of kernel equations, which is the main technical challenge, is shown there by an improved successive approximation method. 



In this paper, based on the results for the linear case \cite{hu-dimeglio-vazquez-krstic}, we will use the linearized feedback control to stabilize the nonlinear system as it is mentioned in \cite{JMC-Vazquez}. Although the target system is a little different from the one in \cite{JMC-Vazquez} with a linear term involved in the equations,  thanks to its special structure, we show  that all the procedures to handle nonlinearities in \cite{JMC-Vazquez} can be also adapted in this paper with more technical developments. Let us recall some definitions and statements  \cite{JMC-Vazquez}. 
Define the norms 
\begin{align*}
&\|u(t,\cdot)\|_{H^1}=\|u(t,\cdot)\|_{L^2}+\|u_x(t,\cdot)\|_{L^2},\\
&\|u(t,\cdot)\|_{H^2}=\|u(t,\cdot)\|_{H^1}+\|u_{xx}(t,\cdot)\|_{L^2}.
\end{align*}
Our main result is given by
\begin{thm}\label{main resultNL}
 Under the assumptions  in $\S$\ref{problem}, suppose furthermore that $C^1$ compatibility conditions   are satisfied at the point $(t, x) = (0, 0)$, there exists a continuous linear feedback control laws $\mathcal{B}: (H^2(0,1))^n\rightarrow \mathbb{R}^m$, satisfying the $C^1$ compatibility conditions at the points $(t,x)=(0,1)$,  then for every $\lambda>0$, there exist $\delta>0$ and $c>0$,  such that the mixed initial-boundary value problem (\ref{3.1}), (\ref{3.10}), (\ref{3.11}) and (\ref{initial}) with $H(t)=\mathcal{B}(u(t,\cdot))$ admits a unique $C^0([0,\infty), (H^2(0,1))^n)$ solution $u=u(t,x)$, which verifies 
\begin{align}\label{mainNL}
\|u(t,\cdot)\|_{H^2}\leq c e^{-\lambda t}\|\phi\|_{H^2},
\end{align}
provided that $\|\phi\|_{H^2}\leq \delta$.
 \end{thm}

\begin{rem}
The $C^1$ compatibility conditions at the point $(t,x)=(0,0)$ are given by
\begin{align}
&\phi_s(0)=G_s(\phi_1(0),\cdots,\phi_m(0))\ &s=m+1,\cdots,n, \label{condition1}\\
&f_s(0,\phi(0))-\sum_{j=1}^n a_{sj}(0,\phi(0))\phi_j'(0)=\notag&\\ 
&\ \ \sum_{r=1}^m\frac{\partial G_s}{\partial u_r}(\phi_1(0),\cdots,\phi_m(0))\cdot\Big(f_r(0,\phi(0))-\sum_{j=1}^n a_{rj}(0,\phi(0))\phi_j'(0)\Big)\ \ \ &s=m+1,\cdots,n. \label{condition2}
\end{align}
The $C^1$ compatibility conditions at the point $(t,x)=(0,1)$ are similar.
\end{rem}

\begin{rem}
For convenience, we always assume that the feedback controls $H(t)=\mathcal{B}(u(t,\cdot))$ satisfy the $C^1$ compatibility conditions at the point $(t,x)=(0,1)$. However, if this property fails, one can add some dynamic terms to the controllers (see also Remark\ \ref{rem artificial condition} and \cite[Section 4]{JMC-Vazquez}). 
\end{rem}

The rest of this paper is organized as follows. In $\S$\ref{review linear case}, we review a former result on the boundary backstepping controls for $n\times n$ linear hyperbolic system. Besides, we design a Lyapunov function to stabilize the linear system in the $L^2$ norm. In $\S$\ref{nonlinear case}, we input the corresponding  linearized closed-loop control to the original nonlinear system and give the feedback control design. In $\S$\ref{proof of main result}, we prove  exponential stability of zero equilibrium for the quasilinear system by using the Control Lyapunov Function method.  We finally  include two appendices with some technical details.


\section{Preliminaries--Linear Case}\label{review linear case}
In this section, we review the results on  stabilization of $n\times n$ hyperbolic linear system by using the backstepping method.  Similar to the situation in \cite{JMC-Vazquez}, this procedure can be applied to locally stabilize the original nonlinear system. Consider the following $ n\times n$  hyperbolic systems
\begin{align}\label{1.1}
w_t(t,x)+\Lambda(x) w_x(t,x)=\Sigma(x) w(t,x),
\end{align}
where, $w=(w_1,\cdots,w_n)^T$ is a vector function of $(t,x)$, $\Lambda$: $[0,1]\rightarrow \mathcal{M}_{n,n}(\mathbb{R})$ is an $n\times n$ $C^2$ diagonal matrix, i.e.
\begin{align}
\Lambda(x)=\left(\begin{array}{cc}\Lambda_-(x) & 0 \\0 & \Lambda_+(x)\end{array}\right),
\end{align} 
in which $\Lambda_-(x):=\mathrm{diag}(\lambda_1(x),\cdots,\lambda_m(x))$ and $\Lambda_+(x):=\mathrm{diag}(\lambda_{m+1}(x),\cdots,\lambda_n(x))$ are  diagonal submatrices, without loss of generality, satisfying 
\begin{align}\label{1.3}
\lambda_1(x)<\cdots<\lambda_m(x)<0< \lambda_{m+1}(x)<\cdots<\lambda_n(x),\forall x\in[0,1].
\end{align}
On the other hand, $\Sigma:[0,1]\rightarrow\mathcal{M}_{n,n}(\mathbb{R})$ is a $n\times n$ matrix  with  
\begin{align}
\Sigma(x)=\left(\begin{array}{cc}\Sigma^{--}(x)& \Sigma^{-+}(x)\\ \Sigma^{+-}(x)& \Sigma^{++}(x)\end{array}\right),
\end{align}
in which $\Sigma^{--}\in \mathcal{M}_{m,m}(\mathbb{R})$, $\Sigma^{-+}\in \mathcal{M}_{m,n-m}(\mathbb{R})$, $\Sigma^{+-}\in \mathcal{M}_{n-m,m}(\mathbb{R})$ and $\Sigma^{--}\in \mathcal{M}_{n-m,n-m}(\mathbb{R})$
are all $C^2$ submatrices with respect to $x$. Moreover,  for any $i=1,\cdots,n$, we assume that 
\begin{align}
\Sigma_{ii}(x)\equiv 0,\ \ \forall x\in[0,1].
\end{align}
The boundary conditions for the linear hyperbolic system (\ref{1.1}) are given by
\begin{align}\label{1.5}
x=0:\ w_+(t,0)=Q w_-(t,0),
\end{align}
and
\begin{align}\label{1.6}
\begin{aligned}
x=1:\ &w_-(t,1)=U(t).
\end{aligned}
\end{align}
where $w_-\in \mathbb{R}^m, w_+\in \mathbb{R}^{n-m}$ are defined by requiring that $w:= (w_-,w_+)^T$, $U=(U_1,\cdots,U_m)^T$ are boundary feedback controls, $Q\in \mathcal{M}_{n-m,m}$ is a constant matrix. 
Our purpose in this section is to find
a full-state feedback control law for U(t) to ensure that the closed-loop system (\ref{1.1}), (\ref{1.5})-(\ref{1.6}) 
is globally asymptotically stable in the $L^2$ norm, which is defined by $\|w(t,\cdot)\|_{L^2}=\sqrt{\sum\limits_{i=1}^n\int_0^1w^2_i(t,x)dx}$.

\subsection{Target System}
In Section \ref{sec3.2}, it will be shown that  we can transform the system (\ref{1.1}), (\ref{1.5})-(\ref{1.6}) into the following cascade system
\begin{align}\label{1.7}
\gamma_t(t,x)+\Lambda(x)\gamma_x(t,x)=G(x)\gamma(t,0)
\end{align}
with the boundary conditions
\begin{align}\label{1.8}
x=0:\ \ \gamma_+(t,0)=Q\gamma_-(t,0)
\end{align}
and
\begin{align}\label{1.9}
\begin{aligned}
x=1:\ &\gamma_-(t,1)=0,
\end{aligned}
\end{align}
where  $\gamma_-\in \mathbb{R}^m, \gamma_+\in \mathbb{R}^{n-m}$ are defined by requiring that $\gamma:= (\gamma_-,\gamma_+)^T$, ~$G$ is a lower triangular matrix with following structure
\begin{align}
G(x)=\left(\begin{array}{cc}\mathcal{G}_1(x) & 0 \\\mathcal{G}_2(x) & 0\end{array}\right),
\end{align}
in which $\mathcal{G}_1\in \mathcal{M}_{m,m}(\mathbb{R})$ is a lower triangular matrix, i.e.
\begin{align}\label{2.12s}
	\mathcal{G}_1(x)&=\begin{pmatrix}
		0&\cdots&\cdots&0\\
		g_{2,1}(x)&\ddots&\ddots&\vdots\\
		\vdots&\ddots&\ddots&\vdots\\
		g_{m,1}(x)&\cdots&g_{m,m-1}(x)&0
	\end{pmatrix},
\end{align}
 and $\mathcal{G}_2(x)\in \mathcal{M}_{n-m,m}(\mathbb{R})$. The coefficients of both $\mathcal{G}_1$ and $\mathcal{G}_2$ are to be determined in $\S$2.2. Next, we prove that the cascade system (\ref{1.7})-(\ref{1.9}) verifies the following proposition.
\begin{proposition}\label{pro2.1}
For any given  matrix function $G(\cdot)\in C^{1}(0,1)$, the mixed initial-boundary value problem (\ref{1.7})-(\ref{1.9})  with initial condition 
\begin{align}
t=0: \gamma(0,x)=\gamma_0(x),
\end{align}
where $\gamma_0\in (L^2(0,1))^n$  admits a $C^0([0,\infty);(L^2(0,1))^n)$ solution $\gamma=\gamma(t,x)$, which is globally exponentially stable in the $L^2$ norm, i.e. for every $\lambda>0$, there exists $c>0$ such that
\begin{align}
\|\gamma(t,\cdot)\|_{L^2}\leq c e^{-\lambda t}\|\gamma_0\|_{L^2}.
\end{align}
In fact, this solution vanishes in finite time $t>t_F$, where $t_F$ is given by
 \begin{align}\label{t_F}
t_F= \int_0^1 \frac{1}{\lambda_{m+1}(s)}+\sum_{r=1}^m\frac{1}{|\lambda_r(s)|}ds.
\end{align}
\end{proposition} 
\textbf{Proof.}
Equations (\ref{1.7}) can be rewritten as
\begin{equation}\label{2.16s}
\begin{aligned}
\partial_t\gamma_-(t,x)+\Lambda_-(x)\partial_x\gamma_-(t,x)=\mathcal{G}_1(x)\gamma_-(t,0),\\
\partial_t\gamma_+(t,x)+\Lambda_+(x)\partial_x\gamma_+(t,x)=\mathcal{G}_2(x)\gamma_-(t,0),
\end{aligned}
\end{equation}  
	then consider the following  Lyapunov functional
\begin{align}
	V_0(t)&=\int_0^1   e^{-\delta x} \gamma_+(t,x)^T\left(\Lambda_+(x)\right)^{-1}\gamma_+(t,x) dx-\int_0^1 e^{\delta x} \gamma_-(t,x)^TB\left(\Lambda_-(x)\right)^{-1}\gamma_-(t,x)dx,  \label{eq:V}
\end{align}
where $\delta>0$ is a parameter, $B=\mathrm{diag}(b_1,\cdots,b_m)$ (with $b_{r}>0, r=1,\cdots,m$) whose coefficients are to be determined. Obviously, $\sqrt{V_0}$ is a norm equivalent to $\|\gamma(t,\cdot)\|_{L^2}$. Differentiating $V_0$ with respect to $t$ and integrating by parts yields
\begin{align}
\notag	\dot{V}_0(t)= I+II+III+IV
\end{align}
with
\begin{align*}
I&=\left[- e^{-\delta x}\gamma_+(t,x)^T\gamma_+(t,x)+e^{\delta x}\gamma_-(t,x)^TB\gamma_-(t,x)\right]_0^1,\\
II&=-\int_0^1 \delta e^{-\delta x}\gamma_+(t,x)^T\gamma_+(t,x)dx-\int_0^1 \delta e^{\delta x}\gamma_-(t,x)^TB\gamma_-(t,x)dx,\\
III&=2\int_0^1 e^{-\delta x}\gamma_+(t,x)^T\left(\Lambda_+(x)\right)^{-1}\mathcal{G}_2(x)\gamma_-(t,0)dx,\\ 
IV&=-2\int_0^1e^{\delta x}\gamma_-(t,x)^TB\left(\Lambda_-(x)\right)^{-1}\mathcal{G}_1(x)\gamma_-(t,0)dx.
\end{align*}
Noting the boundary conditions (\ref{1.8})-(\ref{1.9}), we have that
\begin{align}
&\begin{aligned}
I &=-e^{\delta}\gamma_+(t,1)^T\gamma_+(t,1)-\gamma_-(t,0)^T\Big(B-Q^TQ\Big)\gamma_-(t,0),
\end{aligned}\\
&\begin{aligned}\label{III}
III&\leq   \int_0^1  e^{-\delta x}\gamma_+(t,x)^T\gamma_+(t,x)dx+\gamma_-(t,0)^T\int_0^1e^{-\delta x}\mathcal{G}^T_2(x)\left(\Lambda_+(x)\right)^{-2}\mathcal{G}_2(x)dx\gamma_-(t,0)\\
&\leq   \int_0^1  e^{-\delta x}\gamma_+(t,x)^T\gamma_+(t,x)dx+\gamma_-(t,0)^T\int_0^1\mathcal{G}^T_2(x)\left(\Lambda_+(x)\right)^{-2}\mathcal{G}_2(x)dx\gamma_-(t,0),\\
\end{aligned}\\
&\begin{aligned}\label{IV}
 IV &=-2 \int_0^1e^{\delta x}\sum_{m\geq i>j\geq 1}\gamma_i(t,x)\frac{b_i}{\Lambda_i(x)}g_{ij}(x)\gamma_j(t,0)dx\\
 &\leq -M\int_0^1 e^{\delta x}\sum_{m\geq i>j\geq 1}\frac{b_i}{\Lambda_i(x)} \gamma_i^2(t,x)dx-M\int_0^1e^{\delta x}\sum_{m\geq i>j\geq 1}\frac{b_i}{\Lambda_i(x)}\gamma_j^2(t,0)dx\\
 &\leq -M\int_0^1 e^{\delta x}\sum_{m\geq i>j\geq 1}\frac{b_i}{\Lambda_i(x)} \gamma_i^2(t,x)dx+M\mu e^{\delta}\gamma_-(t,0)^T\mathcal{C}\gamma_-(t,0)\\
 &\leq -m M\int_0^1 e^{\delta x}\sum_{i=2}^m\frac{b_i}{\Lambda_i(x)} \gamma_i^2(t,x)dx+M\mu e^{\delta}\gamma_-(t,0)^T\mathcal{C}\gamma_-(t,0)\\
 &\leq-mM\int_0^1e^{\delta x}\gamma_-(t,x)^TB\left(\Lambda_-(x)\right)^{-1}\gamma_-(t,x)dx+M\mu e^{\delta}\gamma_-(t,0)^T\mathcal{C}\gamma_-(t,0),
\end{aligned}
\end{align}
in which
\begin{align}\label{mathcalC}
M:=\|G\|_{L^{\infty}},\ \ \mathcal{C}:=\mathrm{diag}(\mathcal{C}_1,\cdots,\mathcal{C}_m)
\end{align}
with
\begin{align}
\mathcal{C}_r:=\begin{cases}
\sum\limits_{j=r+1}^mb_j,&1\leq r\leq m-1\\
0, & r=m,
\end{cases}
\end{align}
 and 
\begin{align}\label{mu}
\mu:=\max_{i} \bigg\{\frac{1}{\|\Lambda_i\|_{C^0}}\bigg\}.
\end{align}

Let
\begin{align}
P=Q^TQ+\int_0^1\mathcal{G}^T_2(x)\left(\Lambda_-(x)\right)^{-2}\mathcal{G}_2(x)dx.
\end{align} 
There exists a diagonal matrix $S=\mathrm{diag} (s_1,\cdots,s_m )$ with $s_r>0(\ r=1,\cdots,m)$ being large enough, such that 
\begin{align}\label{S}
P\prec S,
\end{align}
where $P\prec S$ denotes that $S-P$ is a positive-definite matrix. This yields
\begin{align*}
\dot{V}_0(t)\leq& -\gamma_-(t,0)^T\Big(B-S-M\mu e^{\delta}\mathcal{C}\Big)\gamma_-(t,0)-(\delta-1) \int_0^1 e^{-\delta x}\gamma_+(t,x)^T\gamma_+(t,x)dx\\
&-(\delta-mM\mu)\int_0^1e^{\delta x}\gamma_-(t,x)^TB\gamma_-(t,x)dx.
\end{align*}

Thus, for any given $\lambda>0$, picking
\begin{align}
	&\delta>\max \left\{\lambda\mu+mM\mu, \lambda\mu+1\right\}\\
	&\label{chooseB}b_r>\begin{cases}
M\mu e^{\delta}\sum\limits_{j=r+1}^mb_j+s_r,&1\leq r\leq m-1\\
s_m, & r=m,
\end{cases} 
\end{align}
we have
\begin{align}\label{3.17op}
\dot{V}_0\leq -\lambda V_0
\end{align}
where $\lambda$ can be chosen as large as desired. It is easy to see that Parameter matrix $B$ does exist, since one can easily check (\ref{chooseB}) by induction. This shows exponential stability
of $\gamma$ system. 

To show finite-time convergence to the origin, one can find the explicit solution of (\ref{1.7})-(\ref{1.9}) as follows. Define
\begin{align}
\phi_i(x)=\int_0^x\frac{1}{|\lambda_i(\xi)|}d\xi,\ \ 1\leq i\leq n.
\end{align}
Notice that every $\phi_i(1\leq i\leq n)$ is monotonically increasing $C^2$ functions of $x$, and thus invertible.
With  the same statement in \cite{JMC-Vazquez} and noting (\ref{1.7})-(\ref{2.12s}), one can express the explicit solution of $\gamma_1$ by 
\begin{align}\label{2.8z}
\begin{aligned}
\gamma_1(t,x)=\begin{cases}
\gamma_{1}(0,\phi_1^{-1}(\phi_1(x)+t)) &\text{if } t< \phi_1(1)-\phi_1(x),\\
0&\text{if } t\geq \phi_1(1)-\phi_1(x).
\end{cases}
\end{aligned}
\end{align}
Notice in particular that~$\gamma_1$ is identically zero for~$t\geq \phi_1(1)$. From (\ref{1.7}) and (\ref{2.12s}), we obtain that $\gamma_2(t,x)$ satisfies the following equation for~$t\geq \phi_1(1)$
\begin{align}\label{2.9z}
\partial_t\gamma_{2}(t,x)+\lambda_2(x)\partial_x\gamma_{2}(t,x)=0,
\end{align}
with
\begin{align}
\gamma_2(t,1)=0,
\end{align}
which ensures the explicit expression of~$\gamma_2(t,x)$ to be
\begin{align}
\begin{aligned}
\gamma_2(t,x)=\begin{cases}
\gamma_{2}(\phi_1(1),\phi_2^{-1}(\phi_2(x)+t)) &\text{if } \phi_1(1)<t< \phi_1(1)+\phi_2(1)-\phi_2(x),\\
0&\text{if } t\geq \phi_1(1)+\phi_2(1)-\phi_2(x).
\end{cases}
\end{aligned}
\end{align}
Therefore, by induction, one has that $\gamma_r(t,x)(2\leq r\leq m)$ satisfies the following equations, for $ t>\sum\limits_{k=1}^{r-1}\phi_k(1)$,
\begin{align}
\partial_t\gamma_{r}(t,x)+\lambda_r(x)\partial_x\gamma_{r}(t,x)=0,
\end{align}
with the boundary condition
\begin{align}
\gamma_r(t,1)=0.
\end{align}
 Thus, when $t> \sum\limits_{k=1}^{r-1}\phi_k(1)$, we have
\begin{align}
\begin{aligned}
\gamma_r(t,x)=\begin{cases}
\gamma_{r}(\sum\limits_{k=1}^{r-1}\phi_k(1),\phi_r^{-1}(\phi_r(x)+t)) &\text{if } \sum\limits_{k=1}^{r-1}\phi_k(1)<t< \sum\limits_{k=1}^{r}\phi_k(1)-\phi_r(x),\\
0&\text{if } t\geq  \sum\limits_{k=1}^{r}\phi_k(1)-\phi_r(x).
\end{cases}
\end{aligned}
\end{align}
This yields that $\gamma_-(t,x)\equiv 0\ \big(t>\sum\limits_{k=1}^{m}\phi_k(1)\big)$. From the time $t=\sum\limits_{k=1}^{m}\phi_k(1)$ on, we find  $\gamma_+$ becomes the solution of the following system
\begin{align}\label{gamma+}
\partial_t\gamma_+(t,x)+\Lambda_+(x)\partial_x\gamma_+(t,x)=0
\end{align}
with
\begin{align}\label{gamma+boundary}
x=0:\gamma_+(t,0)\equiv 0.
\end{align}
Since (\ref{gamma+})-(\ref{gamma+boundary}) is a completely decoupled system, by the characteristic method,  after $t=t_F$, where 
\begin{align}
t_F=\phi_{m+1}(1)+\sum_{r=1}^m\phi_r(1)=\int_0^1 \frac{1}{\lambda_{m+1}(s)}+\sum_{r=1}^m\frac{1}{|\lambda_r(s)|}ds,
\end{align}
one can see that $\gamma_+(t,x)\equiv 0(t\geq t_F)$, which concludes the Proof of Proposition \ref{pro2.1}.\cqfd


\subsection{Backstepping transformation and Kernel Equations}\label{sec3.2}
To map the original system (\ref{1.1}) into the target system (\ref{1.7}), we use the following Volterra transformation of the second kind,
which is similar to the one in \cite{JMC-Vazquez} and \cite{Meglio}:
\begin{align}\label{2.1}
\gamma(t,x)=w(t,x)-\int_0^x K(x,\xi) w(t,\xi) d\xi.
\end{align}
We point out here that this transformation yields that $w(t,0)\equiv \gamma(t,0)\ (\forall t>0)$, which is crucial to design our feedback law.

Utilizing (\ref{1.1}) and  straightforward computations, one can show that
\begin{align}
\begin{aligned}
\gamma_t+\Lambda(x)\gamma_x=&-\int_0^x \big(K_{\xi}(x,\xi)\Lambda(\xi)+\Lambda(x) K_x(x,\xi)+K(x,\xi)\Sigma(\xi)+K(x,\xi)\Lambda_{\xi}(\xi)\big) w(t,\xi) d\xi\\
&+\big(\Sigma(x)+K(x,x)\Lambda(x)-\Lambda(x) K(x,x)\big) w(t,x)-K(x,0)\Lambda(0)\left(\begin{array}{cc}I&0\\Q&0\end{array}\right) w(t,0).\\
\end{aligned}
\end{align}
The original system (\ref{1.1}) is mapped into the target system (\ref{1.7}) if one has the following  kernel equations: 
\begin{align}
&\Lambda(x) K_x(x,\xi)+K_{\xi}(x,\xi)\Lambda(\xi)+K(x,\xi)\Sigma(\xi)+K(x,\xi)\Lambda_{\xi}(\xi)=0\label{2.6}\\
&\Sigma(x)+K(x,x)\Lambda(x)-\Lambda(x) K(x,x)=0\label{boundaryKxx}\\
&G(x)=-K(x,0)\Lambda(0)\left(\begin{array}{cc}I&0\\Q&0\end{array}\right)\label{2.29}
\end{align}
Developing equations~\eqref{2.6}--\eqref{2.29} leads to the following set of kernel PDEs
\begin{align}
	\lambda_i(x) \partial_x K_{ij}(x,\xi)+\lambda_j(\xi) \partial_\xi K_{ij}(x,\xi) &=- \sum\limits_{k=1}^n\big(\sigma_{kj}(\xi)+\delta_{kj}\lambda_j'(\xi)\big)K_{ik}(x,\xi)\label{eq:developedKernelK}
	\end{align}
along with the following set of boundary conditions
\begin{align}
	K_{ij} (x,x)&=  \cfrac{\sigma_{ij}(x)}{\lambda_i(x)-\lambda_j(x)} \stackrel{\Delta}{=}k_{ij}(x) &\text{for } 1\leq i, j\leq n(i\ne j), \label{eq:hypotenuseK}\\
		 K_{ij}(x,0) &=-\frac{1}{\lambda_j(0)}\sum\limits_{k=1}^{n-m} \lambda_{m+k}(0) K_{i,m+k}(x,0)q_{k,j} & \text{for } 1 \leq i \leq j \leq m.\label{eq:modifx0boundary}
\end{align}
To ensure well-posedness of the kernel equations, 
we add the following artificial boundary conditions for $K_{ij}(m\geq i>j\geq 1, n\geq j>i\geq m+1)$ on $x=1$:
\begin{align}\label{eq:artificialboundary1}
K_{ij}(1,\xi)=k^{(1)}_{ij}(\xi), \ \text{for} \ \ 1&\leq j < i \leq m \ \cup \ m+1\leq i<j\leq n,
\end{align}
and the boundary conditions for $K_{ij}(n\geq i\geq j\geq m+1)$ on $\xi=0$:
\begin{align}\label{eq:artificialboundary}
K_{ij}(x,0)=k_{ij}^{(2)}(x), \ \text{for} \ \ m+1&\leq j\leq i \leq n.
\end{align}
where $k_{ij}^{(1)}$ and $k_{ij}^{(2)}$ are chosen as functions of $C^{\infty}[0,1]$ satisfying the $C^1$ compatibility conditions at the point $(x,\xi)=(1,1)$ (see Remark \ref{rem_compa}).
 The equations evolve in the triangular domain $\mathcal{T}=\{(x,\xi):0\leq \xi\leq x\leq 1\}.$ By Theorem \ref{theorem appendix}, one finds that there exists a unique piecewise $C^{2}(\mathcal{T})$ solution $K(x,\xi)$ to (\ref{eq:developedKernelK})-(\ref{eq:artificialboundary}) with $K(x,x),\ K(x,0) \in C^{1}(0,1)$, provided that $\sigma_{ij}(x)$ are $C^2[0,1]$, $\lambda_i(x)$ are $C^2[0,1]$.  While $G(x)\in C^1$ (with bounded $C^1$ norm) is  given by (\ref{2.29}) under the well-posedness of $K(x,0)$, which is proved in Theorem \ref{theorem appendix}.
 
 \begin{rem}\label{rem_compa}
 The $C^1$ compatibility conditions at the point $(x,\xi)=(1,1)$ are given by
 \begin{align}
&k_{ij}(1)=k_{ij}^{(1)}(1), \ \text{for} \ \ 1\leq j < i \leq m \ \cup \ m+1\leq i<j\leq n,\label{designcomconC0}\\
&\dot{k}^{(1)}_{ij}(1)=\frac{\lambda_i(1)k'_{ij}(1)+\sum\limits_{k=1}^n\big(\sigma_{kj}(1)+\delta_{kj}\lambda_j'(1)\big)k_{ik}(1)}{\lambda_i(1)-\lambda_j(1)},\ \text{for}\ 1\leq j < i \leq m \cup m+1\leq i<j\leq n.\label{designcomconC1}
\end{align}
 \end{rem}


\subsection{The inverse transformation and stabilization for linear system}
Transformation (\ref{2.1}) is a classical Volterra equation of the second kind, one can check from Theorem \ref{theorem appendix2} that there exists a unique piecewise $C^2(\mathcal{T})$ matrix function $L(x,\xi)$ such that
\begin{align}\label{4.1}
w(t,x)=\gamma(t,x)+\int_0^xL(x,\xi)\gamma(t,\xi)d\xi.
\end{align}
From the transformation (\ref{2.1}) evaluated at $x=1$, one gets the following feedback control laws
\begin{align}\label{4.2}
U_i(t)=\int_0^1 \sum_{j=1}^n K_{ij}(1,\xi) w_j(t,\xi) d\xi,\ \ (i=1,\cdots,m),
\end{align}
which immediately leads to our feedback stabilization result for the linear system as follows:
 \begin{thm}
The mixed initial-boundary value problem (\ref{1.1}) with the boundary conditions (\ref{1.5}), the feedback control law (\ref{4.2}) and initial condition
\begin{align}
t=0: w(0,x)=w_0(x),
\end{align}
in which $w_0\in (L^2(0,1))^n$, admits a $(L^2(0,1))^n$ solution $w=w(t,x)$. Moreover, for every $\eta>0$, there exists $c>0$ such that
\begin{align}\label{main}
\|w(\cdot,t)\|_{L^2}\leq c e^{-\eta t}\|w_0\|_{L^2}.
\end{align}
In fact, $w$ vanishes in finite time $t>t_F$, where $t_F$ is given by (\ref{t_F}).
\end{thm}

\begin{rem}
If we focus on the linear problem, $\Lambda$ and $\Sigma$ can be assumed to be $C^1([0,1])$ and $C^0([0,1])$ functions. The corresponding kernels $K$ and $L$ are then both functions of $L^{\infty}(\mathcal{T})$.
\end{rem}

\section{Backstepping boundary control design for nonlinear system}\label{nonlinear case}
As mentioned in \cite{JMC-Vazquez}, we wish  the linear controller (\ref{4.2}) designed by backstepping method to work locally for the corresponding nonlinear system. Let us show that this is indeed the case. Introduce 
\begin{align}
&\varphi_i(x):=\exp{\Big(-\int_{0}^x\frac{f_{ii}(s)}{\Lambda_{i}(s)}ds\Big)}\ \ \ i=1,\cdots,n.
\end{align}
One can make the following coordinates transformation 
\begin{equation}
{w}(t,x)=\left(\begin{array}{ccc}\varphi_1(x)&  &  \\ & \ddots &  \\ &  & \varphi_n(x)\end{array}\right)u(t,x)=\Phi(x) u(t,x).
\end{equation}
Then the original control system $u$ is transformed into the following system expressed in the new coordinates:
\begin{align}\label{new}
w_t(t,x)+\overline{A}(x,w) {w}_x(t,x)=\widetilde{F}(x,w),
\end{align}
in which 
\begin{align}
&\overline{A}(x,w) =\Phi(x)A(x,\Phi^{-1}(x)w)\Phi^{-1}(x),\\
&\widetilde{F}(x,w)=\Phi(x)F(x,\Phi^{-1}(x)w)-\overline{A}(x,w)\left(\begin{array}{ccc}\frac{f_{11}(x)}{\Lambda_1(x)}&  &  \\ & \ddots &  \\ &  & \frac{f_{nn}(x)}{\Lambda_n(x)}\end{array}\right)w.
\end{align}
Obviously, one can check that
\begin{align}
&\widetilde{F}(x,0)=0,\\
& \overline{A}(x,0)=\Phi(x)A(x,0)\Phi^{-1}(x)=A(x,0).
\end{align} 
Moreover, define
\begin{align}
\Sigma(x)=\frac{\partial \widetilde{F}(x,w)}{\partial w}\bigg|_{w=0},
\end{align}
we have that
\begin{align}\label{SigmaNL}
\Sigma_{ij}(x)=\begin{cases}
\frac{\varphi_i(x)}{\varphi_j(x)}f_{ij}(x),i\ne j,\\
0, i=j.
\end{cases}
\end{align}
Therefore, we may rewrite (\ref{new}) as a linear system with the same structure as (\ref{1.1}) plus nonlinear terms:
\begin{align}\label{final}
w_t(t,x)+{\Lambda}(x) {w}_x(t,x)=\Sigma(x)w(t,x)+\Lambda_{NL}(x,w)w_x(t,x)+f_{NL}(x,w),
\end{align}
where
\begin{align}\label{LambdaNL}
\Lambda(x)=A(x,0),
\end{align}
and
\begin{align}
\Lambda_{NL}(x,w)=\Lambda(x)-\overline{A}(x,w),\  f_{NL}(x,w)=\widetilde{F}(x,w)-\Sigma(x)w(t,x).
\end{align}
For the boundary conditions of the system (\ref{final}), defining 
\begin{align}
Q=\Big(\frac{\partial G_s}{\partial u_r}\Big)_{(n-m)\times m}\bigg|_{u=0}\ \ \text{and} \ G_{NL}(w_-(t,0))=G(w_-(t,0))-Qw_-(t,0),
\end{align}
one obtains that
\begin{align}\label{finalbc}
x=0:w_+(t,0)=Qw_-(t,0)+G_{NL}(w_-(t,0))
\end{align} 
and
\begin{align}\label{finalbc1}
x=1:w_-(t,1)=U(t),
\end{align} 
where
\begin{align}
U(t)=\left(\begin{array}{ccc}\varphi_1(1)&  &  \\ & \ddots &  \\ &  & \varphi_m(1)\end{array}\right)H(t)=\widetilde{\Phi}(1) H(t).
\end{align}

It is easily verified that 
\begin{align}
\Lambda(x,0)=0,\ \ f_{NL}(x,0)=\frac{\partial f_{NL}}{\partial w}(x,0)=0
\end{align}
and
\begin{align}
G_{NL}(0)=\frac{\partial G_{NL}}{\partial w}(0)=0.
\end{align}
Thus, the feedback control law  can be chosen as
\begin{align}\label{final control}
h_r(t)=\widetilde{\Phi}_{rr}^{-1}(1) U_r(t)=\widetilde{\Phi}_{rr}^{-1}(1) \int_0^1 \sum_{j=1}^nK_{rj}(1,\xi) \widetilde{\Phi}_{jj}(\xi)u_j(t,\xi) d\xi, \ \ r=1,\cdots,m,
\end{align}
where the kernels are computed from (\ref{eq:developedKernelK})--(\ref{eq:artificialboundary}) with the coefficients $\Sigma(x)$ and $\Lambda(x)$ obtained from (\ref{SigmaNL}) and (\ref{LambdaNL}). One easily verifies that under the assumptions of $\S$\ref{problem}, both $\Sigma$ and $\Lambda$ are functions of $C^2$.

\begin{rem}\label{rem artificial condition}
The $C^1$ compatibility conditions at the point $(t,x)=(0,1)$ for system (\ref{3.1}) with boundary conditions (\ref{finalbc1}) should be
\begin{align}
&\phi_r(1)= \sum_{j=1}^n\int_0^1 \tilde{k}_{rj}(\xi) \phi_j(\xi) d\xi,\ \ \ r=1,\cdots,m, \label{artificial condition1} \\
&f_r(1,\phi(1))-\sum_{j=1}^n a_{rj}(1,\phi(1))\phi_j'(1)=\notag\\
&\ \ \ \sum_{k=1}^n\int_0^1\tilde{k}_{rk}(\xi)\Big(f_k(1,\phi(1))-\sum_{j=1}^n a_{kj}(1,\phi(1))\phi_j'(1)\Big),\ \ \ r=1,\cdots,m, \label{artificial condition2}
\end{align}
where $\tilde{k}_{rk}(\xi)$ are the elements of the matrix $\widetilde{K}(\xi)$ with
\begin{align}
\widetilde{K}(\xi)=\widetilde{\Phi}^{-1}(1) K(1,\xi)  \widetilde{\Phi}(\xi).
\end{align}

Notice that (\ref{artificial condition1})-(\ref{artificial condition2}) depend on the feedback control design, however, there are no physical reasons that the initial data should satisfy them.  In order to guarantee the initial conditions independent of these artificial conditions,  we, following   \cite{JMC-Vazquez}, modify the boundary controls on $x=1$ as
\begin{align}\label{control design}
x=1: \ \ u_r=h_r(t)+a_r(t)+b_r(t),\ \ \ r=1,\cdots,m,
\end{align}
where $a_r$ and $b_r$ are the state of the following dynamic systems
\begin{align}
\dot{a}_r(t)=-d_r a_r(t),\ \ \dot{b}_r(t)=-\tilde{d}_r b_r(t), \ \ r=1,\cdots,m
\end{align}
with $d_r>0,\tilde{d}_r>0$ and $d_r\ne \tilde{d}_r, r=1,\cdots,m$.  By the modified control designs (\ref{control design}), the compatibility conditions on $x=1$ are rewritten by
\begin{align}
&\phi_r(1)= \sum_{j=1}^n\int_0^1 \tilde{k}_{rj}(\xi) \phi_j(\xi) d\xi+a_r(0)+b_r(0),\ \ \ r=1,\cdots,m, \\
&f_r(1,\phi(1))-\sum_{j=1}^n a_{rj}(1,\phi(1))\phi'_j(1)=\notag\\
&\ \ \ \sum_{k=1}^n\int_0^1\tilde{k}_{rk}(\xi)\Big(f_k(1,\phi(1))-\sum_{j=1}^n a_{kj}(1,\phi(1))\phi'_j(1)\Big)-d_r a_r(0)-\tilde{d}_r b_r(0),\ \  r=1,\cdots,m. 
\end{align}
For any $1\leq r\leq m$, call 
\begin{align}
&\mathcal{P}_r(\phi)=\phi_r(1)- \sum_{j=1}^n\int_0^1 \tilde{k}_{rj}(\xi) \phi_j(\xi) d\xi\\
&\mathcal{M}_r(\phi)=f_r(1,\phi(1))-\sum_{j=1}^n a_{rj}(1,\phi(1))\phi'_j(1)-\sum_{k=1}^n\int_0^1\tilde{k}_{rk}(\xi)\Big(f_k(1,\phi(1))-\sum_{j=1}^n a_{kj}(1,\phi(1))\phi'_j(1)\Big)
\end{align}
Picking
\begin{align}\label{ab}
a_r(0)=-\frac{\mathcal{M}_r(\phi)+\tilde{d}_r\mathcal{P}_r(\phi)}{d_r-\tilde{d}_r},\ \  b_r(0)=\frac{d_r\mathcal{P}_r(\phi)+\mathcal{M}_r(\phi)}{d_r-\tilde{d}_r},
\end{align}
the compatibility conditions are automatically verified. Similar stabilization results as Theorem \ref{main resultNL} are still valid for the closed--loop system (\ref{3.1}), (\ref{3.10}) and (\ref{control design}) (see \cite[Theorem 4.1]{JMC-Vazquez}). In fact, this dynamic extension is designed to avoid restriction for artificial boundary conditions due to the compatibility conditions at the points $(t,x)=(0,1)$, and it has been introduced in \cite{JMC-1999} to deal with the stabilization of the Euler equations of incompressible fluids (see also \cite{VTC-2008}).
\end{rem}


 \section{Proof of Theorem \ref{main resultNL}}\label{proof of main result}
 In this section, we will prove the exponential stability for the  system (\ref{3.1}), (\ref{3.10}) and (\ref{3.11}) under the boundary feedback controls (\ref{final control}) by Control Lyapunov Function method.  The whole proof is divided into the following steps.
 \subsection{Definitions}
 We first define some notations (omitting the time argument):
 \begin{align}
 \|\gamma\|_{\infty}:=\operatorname*{ess\,sup}_{x\in[0, 1]}|\gamma(x)|,\ \ \|\gamma\|_{L^p}:=\Big(\int_0^1|\gamma(\xi)|^pd\xi\Big)^{\frac{1}{p}},\ 1\leq p<+\infty.
 \end{align}
 For a $n\times n$ matrix, denote
 \begin{align}
 |M|:=\max\{\|M\gamma\|_{L^{\infty}}:\gamma\in\mathbb{R}^n, |\gamma|=1\}.
 \end{align}
 For a piecewise kernel matrix $K(x,\xi)$, which is a continuous function on each domain $D_i(i=1,\cdots,\mathcal{S})$, respectively, with 
 \begin{align}
 &\mathcal{T}=\bigcup_{i=1}^\mathcal{S}D_i,\\
 &D_i\cap D_j=\emptyset, (i\ne j).
 \end{align} 
 Let
 \begin{align}
 \|K\|_{\infty}:=\max_{i} \operatorname*{sup}_{(x,\xi)\in D_i}|K(x,\xi)|.
 \end{align}
 As before, we recall the following symbols of  \cite{JMC-Vazquez} for simplicity:
 \begin{align}
 \mathcal{K}[\gamma](t,x)&=\gamma(t,x)-\int_0^xK(x,\xi)\gamma(t,\xi)d\xi,\\
  \mathcal{L}[\gamma](t,x)&=\gamma(t,x)+\int_0^xL(x,\xi)\gamma(t,\xi)d\xi,\\
   \mathcal{K}_1[\gamma](t,x)&=-K(x,x)\gamma(t,x)+\int_0^{x}K_{\xi}(x,\xi)\gamma(t,\xi)d\xi,\\
   \mathcal{K}_2[\gamma](t,x)&=-K(x,x)\gamma(t,x)-\int_0^{x}K_{x}(x,\xi)\gamma(t,\xi)d\xi,\\
   \mathcal{L}_1[\gamma](t,x)&=L(x,x)\gamma(t,x)+\int_0^{x}L_{x}(x,\xi)\gamma(t,\xi)d\xi.
 \end{align}
 
 Define $F_1[\gamma]$ and $F_2[\gamma]$ as
 \begin{align}
 F_1[\gamma]:=\Lambda_{NL}(x,\mathcal{L}[\gamma]),\ \ F_2[\gamma]:=f_{NL}(x,\mathcal{L}[\gamma]).
 \end{align}
To prove our result, we notice that if we apply the (inverse) backstepping transformation (\ref{2.1})
to the nonlinear system (\ref{final}), we obtain the following transformed system
\begin{align}\label{analysis gamma}
\begin{aligned}
&\gamma_t(t,x)+\Lambda(x)\gamma_x(t,x)-G(x)\gamma(t,0)\\
=&\mathcal{K}[\Lambda_{NL}(x,w)w_x]+\mathcal{K}[f_{NL}(x,w)]\\
=&\mathcal{K}[\Lambda_{NL}(x,w)\gamma_x]+\mathcal{K}[\Lambda_{NL}(x,w)\mathcal{L}_1[\gamma]]+\mathcal{K}[f_{NL}(x,w)]\\
=&F_3[\gamma,\gamma_x]+F_4[\gamma],
\end{aligned}
\end{align}
where
\begin{align*}
&F_3=\mathcal{K}[F_1[\gamma]\gamma_x],\\
&F_4=\mathcal{K}[F_1[\gamma]\mathcal{L}_1[\gamma]+F_2[\gamma]].
\end{align*}
The boundary conditions are
\begin{align}
x=0:\gamma_+(t,0)=Q\gamma_-(t,0)+G_{NL}(\gamma_-(t,0))
\end{align} 
and
\begin{align}\label{BCgamma}
x=1:\gamma_-(t,1)=0.
\end{align}

 
Notice that here  we may lose the regularity on the point $(0,0)$ for the kernels $K$ and $L$, which leads both of them  to be discontinuous (see \cite{hu-dimeglio-vazquez-krstic}). However, by the assumptions on the coefficients and applying Theorem \ref{theorem appendix} and Theorem \ref{theorem appendix2}, the direct and inverse transformations (\ref{2.1}) and (\ref{4.1}) have $C^{2}$ piecewise kernels functions. 
Fortunately, differentiating twice with respect to $x$ in these transformations, by the similar argument in \cite{JMC-Vazquez} and \cite[Proposition 3.1]{VTC-2008} as well as the additive property of the integral, it can be shown that the $H^2$ norm of $\gamma$ is equivalent to the $H^2$ norm of $w$. Thus, if we show $H^2$ local stability of the origin for (\ref{analysis gamma})-(\ref{BCgamma}), the same holds for $w$ i.e. $u$.
 
 In order to get the desired $H^2$ estimation  for $\gamma$, the things left are just estimating the growth of $\|\gamma\|_{L^2}$, $\|\gamma_t\|_{L^2}$ and $\|\gamma_{tt}\|_{L^2}$, respectively.
 
 \subsection{Analyzing the growth of $\|\gamma\|_{L^2}$}
 Let
 \begin{align}
 F_3[\gamma,\gamma_x]=\big(F_3^-[\gamma,\gamma_x],F_3^+[\gamma,\gamma_x]\big)^T,\ \ \ \ 
 F_4[\gamma]=\big(F_4^-[\gamma],F_4^+[\gamma]\big)^T.
 \end{align}
 where $F_3^-$ and $F_4^- \in \mathbb{R}^m$, $F_3^+$ and $F_4^+\in \mathbb{R}^{n-m}$.
 
 Define
 \begin{align}
	V_1(t)&=\int_0^1   e^{-\delta x} \gamma_+(t,x)^T\left(\Lambda_+(x)\right)^{-1}\gamma_+(t,x) dx-\int_0^1 e^{\delta x} \gamma_-(t,x)^TB\left(\Lambda_-(x)\right)^{-1}\gamma_-(t,x)dx.  \label{eq:V1}
\end{align}
Differentiating $V_1$ with respect to time and integrating by parts yields
\begin{align}
\notag	\dot{V}_1(t)=V+VI+VII+VIII+IX+X
\end{align}
with
\begin{align*}
V=&\left[- e^{-\delta x}\gamma_+(t,x)^T\gamma_+(t,x)+e^{\delta x}\gamma_-(t,x)^TB\gamma_-(t,x)\right]_0^1,\\
VI=&-\int_0^1 \delta e^{-\delta x}\gamma_+(t,x)^T\gamma_+(t,x)dx-\int_0^1 \delta e^{\delta x}\gamma_-(t,x)^TB\gamma_-(t,x)dx,\\
VII=&2\int_0^1 e^{-\delta x}\gamma_+(t,x)^T\left(\Lambda_+(x)\right)^{-1}\mathcal{G}_2(x)\gamma_-(t,0)dx,\\
VIII=&-2\int_0^1e^{\delta x}\gamma_-(t,x)^TB\left(\Lambda_-(x)\right)^{-1}\mathcal{G}_1(x)\gamma_-(t,0)dx,\\
IX=&2\int_0^1 e^{-\delta x}\gamma_+(t,x)^T\left(\Lambda_+(x)\right)^{-1}\big(F_3^+[\gamma,\gamma_x]+F_4^+[\gamma]\big)dx,\\
X=&-2\int_0^1e^{\delta x}\gamma_-(t,x)^TB\left(\Lambda_-(x)\right)^{-1}\big(F_3^-[\gamma,\gamma_x]+F_4^-[\gamma]\big)dx.
\end{align*}
 By the same argument in \cite {JMC-Vazquez} and noting Lemma \ref{B.2lemma}, we have
\begin{align}
\begin{aligned}
IX+X&\leq K_1\int_0^1|\gamma|(|F_3[\gamma,\gamma_x]|+|F_4[\gamma]|)dx\\
&\leq K_2(\|\gamma_x\|_{\infty}V_1+V_1^{\frac{3}{2}}).
\end{aligned}
\end{align}
Moreover, for $\|\gamma\|_{\infty}\leq \delta$, $|G_{NL}(\gamma_-(t,0))|\leq K_3|\gamma_-(t,0)|$, then 
\begin{align}
\begin{aligned}
V&=- e^{-\delta }\gamma_+(t,1)^T\gamma_+(t,1)+e^{\delta }\gamma_-(t,1)^TB\gamma_-(t,1)+\gamma_+(t,0)^T\gamma_+(t,0)-\gamma_-(t,0)^TB\gamma_-(t,0)\\
&\leq-\gamma_-(t,0)^T\Big(B-Q^TQ-K_3^2I_m\Big)\gamma_-(t,0). 
\end{aligned}
\end{align}
By (\ref{III}) and (\ref{IV}),  one immediately obtains 
\begin{align*}
\dot{V}_1(t)\leq& -\gamma_-(t,0)^T\Big(B-\widetilde{S}-M\mu e^{\delta}\mathcal{C}\Big)\gamma_-(t,0)-(\delta-1) \int_0^1 e^{-\delta x}\gamma_+(t,x)^T\gamma_+(t,x)dx\\
&-(\delta-mM\mu)\int_0^1e^{\delta x}\gamma_-(t,x)^TB\gamma_-(t,x)dx+K_2\big(V_1^{\frac{3}{2}}+\|\gamma_x\|_{\infty}V_1\big), 
\end{align*}
where $M,\ \mathcal{C},\mu$ are given by (\ref{mathcalC}) and (\ref{mu}), $\widetilde{S}:=S+K_3^2I_m$ with $S$ stated in (\ref{S}). Thus, for any given $\lambda_1>0$, picking
\begin{align}
	&\delta>\max \left\{\lambda_1\mu+mM\mu, \lambda_1\mu+1\right\},\\
	&\label{B}b_r:=\begin{cases}
M\mu e^{\delta}\sum\limits_{j=r+1}^mb_j+\tilde{s}_{r},&1\leq r\leq m-1\\
\tilde{s}_m, & r=m,
\end{cases} 
\end{align}
we have the following 
\begin{proposition}\label{proV1}
For any given $\lambda_1>0$, there exists $\delta_1>0$ and $K_2>0$, such that 
\begin{align}
\dot{V}_1\leq -\lambda_1 V_1+K_2\big(V_1^{\frac{3}{2}}+\|\gamma_x\|_{\infty}V_1\big),
\end{align}
provided $\|\gamma\|_{\infty}\leq \delta_1$.
\end{proposition}
 
 \subsection{Analyzing the growth of $\|\gamma_t\|_{L^2}$}
Let $\zeta=\gamma_t$. Taking the partial derivative with $t$ in (\ref{analysis gamma}) yields:
\begin{align}\label{eq:zeta}
\begin{aligned}
\zeta_t(t,x)+(\Lambda(x)-F_1[\gamma])\zeta_x(t,x)-G(x)\zeta(t,0)
=F_5[\gamma,\gamma_x,\zeta]+F_6[\gamma,\zeta],
\end{aligned}
\end{align}
where
\begin{align}
&\begin{aligned}
F_5=\mathcal{K}_1[F_1[\gamma]\zeta]+\int_0^xK(x,\xi)F_{12}[\gamma,\gamma_x]\zeta(\xi)d\xi+K(x,0)\Lambda_{NL}(0,\gamma(0))\zeta(0)+\mathcal{K}[F_{11}[\gamma,\zeta]\gamma_x],
\end{aligned}\\
&\begin{aligned}
F_6=\mathcal{K}[F_{11}[\gamma,\zeta]\mathcal{L}_1[\gamma]]+\mathcal{K}[F_1[\gamma]\mathcal{L}_1[\zeta]]+\mathcal{K}[F_{21}[\gamma,\zeta]],
\end{aligned}
\end{align}
with
\begin{align}
\begin{aligned}
&F_{11}=\frac{\partial \Lambda_{NL}}{\partial \gamma}(x,\mathcal{L}[\gamma])\mathcal{L}[\zeta],\\
&F_{12}=\frac{\partial \Lambda_{NL}}{\partial \gamma}(x,\mathcal{L}[\gamma])(\gamma_x+\mathcal{L}_1[\gamma])+\frac{\partial \Lambda_{NL}}{\partial \gamma}(x,\mathcal{L}[\gamma]),\\
&F_{21}=\frac{\partial f_{NL}}{\partial \gamma}(x,\mathcal{L}[\gamma])\mathcal{L}[\zeta].
\end{aligned}
\end{align}
The boundary conditions are given by
\begin{align}\label{BCzeta0}
x=0:\zeta_+(t,0)=Q\zeta_-(t,0)+\frac{\partial G_{NL}}{\partial \gamma_-}\Big(\gamma_-(t,0)\Big)\zeta_-(t,0)
\end{align} 
and
\begin{align}\label{BCzeta1}
x=1:\zeta_-(t,1)=0,
\end{align} 
in which $\zeta_-\in \mathbb{R}^m, \zeta_+\in \mathbb{R}^{n-m}$ are defined by requiring that $\zeta:= (\zeta_-,\zeta_+)^T$. 

Similarly as in \cite{JMC-Vazquez}, we need the following lemma in order to find a Lyapunov function for $\zeta(t,x)$:
\begin{lemma}
There exists $\delta>0$ such that, for any $\|\gamma\|_{\infty}\leq\delta$, there exists a symmetric matrix $R[\gamma]$ satisfying the identity
\begin{align}\label{symmetry}
R[\gamma](\Lambda(x)-F_1[\gamma])-(\Lambda(x)-F_1[\gamma])^TR[\gamma]=0.
\end{align} 
Moreover, we have that
\begin{align}
&|R[\gamma](x)|\leq c_1+c_2\|\gamma\|_{{\infty}},\\
&\Big|\big((R[\gamma]-D(x))\Lambda(x)\big)_x\Big|\leq c_2(\|\gamma\|_{{\infty}}+\|\gamma_x\|_{{\infty}}),\\
&|(R[\gamma])_t|\leq c_3(|\zeta|+\|\zeta\|_{L^1}),
\end{align}
where $c_1$, $c_2$ and $c_3$ are positive constants, and 
\begin{align}
D(x)=\left(\begin{array}{cc} -e^{\delta x}B(\Lambda_-(x))^{-1}& 0 \\0 & e^{-\delta x}(\Lambda_+(x))^{-1}\end{array}\right).
\end{align}
\end{lemma}
\textbf{Proof:}   Denote $\mathcal{D}_n(x)$ as the set of $n\times n$  diagonal matrices with $C^1$ elements.  Let $\Lambda(x):= \mathrm{diag} (\Lambda_1(x),\cdots,\Lambda_n(x)) \in \mathcal{D}_n(x)$ be such that $\Lambda_i(x)\ne\Lambda_j(x)(i\ne j\forall x\in[0,1])$ holds. Notice that $D\in\mathcal{D}_n(x)$. Based on the proof in \cite[Lemma 4.1]{JMC-Bastin-Andera-2008}, one can easily see that there exist a positive real number $\eta$ and a map $\mathcal{N} : \{M\in \mathcal{M}_{n,n}(\mathbb{R};x); \|M(x)-\Lambda(x)\|_{C^1}<\eta\}\rightarrow \mathcal{S}_n$ of class $C^{\infty}$ such that
\begin{align}
\mathcal{N}(\Lambda(x)) =D(x) ,
\end{align}
and
\begin{align}
\mathcal{N}(M)M-M^T\mathcal{N}(M)=0\ \ \forall M\in \mathcal{M}_{n,n}(\mathbb{R};x),\ \ \|M(x)-\Lambda(x)\|_{C^1}<\eta.
\end{align}
It then suffices to define $R[\gamma]$ by
\begin{align}
R[\gamma]=\mathcal{N}(\Lambda(x)-F_1[\gamma]).
\end{align}
Moreover, by the regularity of $\mathcal{N}$ and Lemma \ref{B.2lemma}--\ref{B.3lemma}, one can show that
\begin{align}
\begin{aligned}
|R[\gamma]|&\leq |D(x)|+|R[\gamma]-D(x)|\\
& \leq c_4+c_5|F_1[\gamma]|\\
&\leq c_4+c_6\|\gamma\|_{{\infty}},
\end{aligned}
\end{align}

\begin{align}
\begin{aligned}
\Big|\big((R[\gamma]-D(x))\Lambda(x)\big)_x\Big|&\leq |(R[\gamma]-D(x))_x\Lambda(x)|+|(R[\gamma]-D(x))\Lambda_x(x)|\\
&\leq c_7|F_{12}|+c_8|F_1|\\
&\leq c_9(\|\gamma\|_{{\infty}}+\|\gamma_x\|_{{\infty}})
\end{aligned}
\end{align}
and
\begin{align}
|R[\gamma]_t|&\leq c_{10}\Big|\frac{\partial F_1[\gamma]}{\partial t}\Big|\\
&\leq c_{10}|F_{11}[\gamma,\zeta]|\\
&\leq c_{11}(|\zeta|+\|\zeta\|_{L^1}).
\end{align}
This concludes the proof of Lemma 4.1.
\cqfd

Define
\begin{align}
V_2(t)&=\int_0^1 \zeta^T(t,x)R[\gamma]\zeta(t,x) dx.    \label{eq:V2}
\end{align}
Using (\ref{symmetry}) and straightforward computations, one can show that
\begin{align*}
\dot{V}_2(t)=XI+XII+XIII+XIV+XV
\end{align*}
with
\begin{align*}
XI=&\int_0^1\zeta^T(t,x)(R[\gamma](\Lambda(x)-F_1[\gamma]))_x\zeta(t,x)dx,\\
XII=&-[\zeta^T(t,x)R[\gamma](\Lambda(x)-F_1[\gamma])\zeta(t,x)]_{x=0}^{x=1},\\
XIII=&\int_0^1\zeta(t,x)(R[\gamma])_t\zeta(t,x)dx,\\
XIV=&2\int_0^1\zeta^T(t,x)R[\gamma]F_5[\gamma,\gamma_x,\zeta,\zeta_x]dx+2\int_0^1\zeta^T(t,x)R[\gamma]F_6[\gamma,\zeta]dx,\\
XV=&2\int_0^1\zeta^T(t,x)R[\gamma]G(x)\zeta(t,0)dx.
\end{align*}
For $XII$ and $XV$, by the boundary conditions (\ref{BCzeta0})--(\ref{BCzeta1}), we have
\begin{align}
\begin{aligned}
XII+XV&=-[\zeta^T(t,x)(D(x)+\Theta[\gamma])(\Lambda(x)-F_1[\gamma])\zeta(t,x)]_{x=0}^{x=1}\\
&\ \ \ \ +2\int_0^1\zeta^T(t,x)(D(x)+\Theta[\gamma])G(x)\zeta(t,0)dx\\
&=-[\zeta^T(t,x)(D(x)\Lambda(x)+\Theta[\gamma]\Lambda(x)-D(x)F_1[\gamma]-\Theta[\gamma]F_1[\gamma])\zeta(t,x)]_{x=0}^{x=1}\\
&\ \ \ \ +2\int_0^1\zeta^T(t,x)D(x)G(x)\zeta(t,0)dx+2\int_0^1\zeta^T(t,x)\Theta[\gamma]G(x)\zeta(t,0)dx\\
&\leq -\zeta_-(t,0)^T\Big(B-\widetilde{S}-M\mu e^{\delta}\mathcal{C}-K_{3}\|\gamma\|_{\infty}I_m\Big)\zeta_-(t,0)\\
&\ \ \ \ \ + \int_0^1 e^{-\delta x}\zeta_+(t,x)^T\zeta_+(t,x)dx+mM\mu\int_0^1e^{\delta x}\zeta_-(t,x)^TB\zeta_-(t,x)dx\\
&\ \ \ \ \  +K_{4}\|\gamma\|_{\infty}V_2.
\end{aligned}
\end{align}
As stated in \cite{JMC-Vazquez}, we obtain
\begin{align}
XI&\leq -\lambda_2V_2 +K_4\|\zeta\|_{L^2}^2(\|\gamma\|_{\infty}+\|\gamma_x\|_{\infty}),\\
XIII&\leq K_{5}\|\zeta\|^2_{L^2}\|\zeta\|_{\infty},\\
XIV&\leq K_{6}\Big(\|\zeta\|^2_{L^2}(\|\gamma\|_{\infty}+\|\gamma_x\|_{\infty})+\|\zeta\|_{L^2}|\zeta(t,0)||\gamma(t,0)|\Big).
\end{align}
Following Lemma \ref{B.5lemma}, we are  in the position to conclude that
\begin{proposition}\label{proV2}
For any given $\lambda_2>0$, there exists $\delta_2>0$ and $K_7>0$,  such that 
\begin{align}
\dot{V}_2\leq -\lambda_2 V_2+K_{7}\big(\|\zeta\|_{\infty}+\|\gamma\|_\infty\big) V_2,
\end{align}
provided that $\|\gamma\|_{\infty}\leq \delta_2$.
\end{proposition}

 \subsection{Analyzing the growth of $\|\gamma_{tt}\|_{L^2}$}
 We next deal with $\|\gamma_{tt}\|_{L^2}$. Define $\theta=\gamma_{tt}$. Taking a partial derivative with respect to $t$ for (\ref{eq:zeta}), one obtains an equation of $\theta$:
 \begin{align}
 \theta_t+[\Lambda(x)-F_1[\gamma]]\theta_x=G(x)\theta(t,0)+F_7[\gamma,\gamma_x,\zeta,\zeta_x,\theta]+F_8[\gamma,\zeta,\theta],
 \end{align}
 where
 \begin{align}
& \begin{aligned}
 F_7&=\mathcal{K}_1[F_{11}[\gamma,\zeta]\zeta]+\int_0^xK(x,\xi)F_{12}[\gamma,\gamma_x]\theta(\xi)d\xi+\mathcal{K}_1[F_1[\gamma]\theta]\\
 &\ \ +\int_0^xK(x,\xi)F_{14}[\gamma,\gamma_x,\zeta,\zeta_x]\zeta(\xi)d\xi+K(x,0)\frac{\partial \Lambda_{NL}}{\partial \gamma}(0,\gamma(0))\zeta(0)\zeta(0)\\
 &\ \ +K(x,0)\Lambda_{NL}(0,\gamma(0))\theta(0)+\mathcal{K}[F_{11}[\gamma,\zeta]\zeta_x]+\mathcal{K}[F_{13}[\gamma,\zeta,\theta]\gamma_x],
 \end{aligned}\\
 &\begin{aligned}
 F_8&=2\mathcal{K}[F_{11}[\gamma,\zeta]\mathcal{L}_1[\zeta]]+\mathcal{K}[F_1[\gamma]\mathcal{L}_1[\theta]]+\mathcal{K}[F_{13}[\gamma,\zeta,\theta]\mathcal{L}_1[\gamma]]+\mathcal{K}[F_{22}[\gamma,\zeta,\theta]]
 \end{aligned}
 \end{align}
 with
 \begin{align}
 &\begin{aligned}
 F_{13}=\frac{\partial \Lambda^2_{NL}}{\partial \gamma^2}(x,\mathcal{L}[\gamma])\mathcal{L}[\zeta]\mathcal{L}[\zeta]+\frac{\partial \Lambda_{NL}}{\partial\gamma}(x,\mathcal{L}[\gamma])\mathcal{L}[\theta],
 \end{aligned}\\
 &\begin{aligned}
 F_{14}=\frac{\partial \Lambda^2_{NL}}{\partial \gamma^2}(x,\mathcal{L}[\gamma])\mathcal{L}[\zeta](\gamma_x+\mathcal{L}_1[\gamma])+\frac{\partial \Lambda_{NL}}{\partial\gamma}(x,\mathcal{L}[\gamma])(\zeta_x+\mathcal{L}_1[\zeta])+\frac{\partial^2 \Lambda_{NL}}{\partial x\partial \gamma}(x,\mathcal{L}[\gamma])\mathcal{L}[\zeta],
 \end{aligned}\\
 &\begin{aligned}
 F_{22}=\frac{\partial^2 f_{NL}}{\partial\gamma^2}(x,\mathcal{L}[\gamma])\mathcal{L}[\zeta]\mathcal{L}[\zeta]+\frac{\partial f_{NL}}{\partial \gamma}(x,\mathcal{L}[\gamma])\mathcal{L}[\theta].
 \end{aligned}
 \end{align}
 The boundary conditions of $\theta$ are given by
\begin{align}\label{BCtheta0}
x=0:\theta_+(t,0)=Q\theta_-(t,0)+\frac{\partial G_{NL}}{\partial \gamma_-}\Big(\gamma_-(t,0)\Big)\theta_-(t,0)+\frac{\partial^2G_{NL}}{\partial \gamma_-^2}\Big(\gamma_-(t,0)\Big)\zeta_-(t,0)\zeta_-(t,0)
\end{align} 
and
\begin{align}\label{BCtheta1}
x=1:\theta_-(t,1)=0.
\end{align} 
where  $\theta_-\in \mathbb{R}^m, \theta_+\in \mathbb{R}^{n-m}$ are defined by requiring that $\theta:= (\theta_-,\theta_+)^T$.

In order to control $\|\theta\|_{L^2}$, we introduce
\begin{align}
V_3(t)&=\int_0^1 \theta^T(t,x)R[\gamma]\theta(t,x) dx,    \label{eq:V3}
\end{align}
then it is easy to see that
\begin{align}\label{V_3prime}
\begin{aligned}
\dot{V}_3(t)=XVI+XVII+XVIII+XIX+XX
\end{aligned}
\end{align}
with
\begin{align*}
XVI=&\int_0^1\theta^T(t,x)(R[\gamma](\Lambda(x)-F_1[\gamma]))_x \theta(t,x)dx,\\
XVII=&-[\theta^T(t,x)R[\gamma](x)(\Lambda(x)-F_1[\gamma](x))\theta(t,x)]_{x=0}^{x=1},\\
XVIII=&+\int_0^1\theta^T(t,x)(R[\gamma])_t\theta(t,x)dx,\\
XIX=&2\int_0^1\theta^T(t,x)R[\gamma]F_7[\gamma,\gamma_x,\zeta,\zeta_x,\theta]dx+2\int_0^1\theta^T(t,x)R[\gamma]F_8[\gamma,\zeta,\theta]dx,\\
XX=&2\int_0^1\theta^T(t,x)R[\gamma]G(x)\theta(t,0)dx.
\end{align*}
Let us first look at the second and the last term of (\ref{V_3prime})(i.e. XVII and XX), by some straight computations, one gets
\begin{align}
\begin{aligned}
XVII+XX\leq& -\theta_-(t,0)^T\Big(B-\widetilde{S}-M\mu e^{\delta}\mathcal{C}-K_{8}\|\gamma\|_{\infty}I_m\Big)\theta_-(t,0)\\
&\ + \int_0^1 e^{-\delta x}\theta_+(t,x)^T\theta_+(t,x)dx+mM\mu\int_0^1e^{\delta x}\theta_-(t,x)^TB\theta_-(t,x)dx\\
&\ +K_{9}\|\gamma\|_{\infty}V_3.\\
\end{aligned}
\end{align}

Then by the same procedures in \cite{JMC-Vazquez}, we have the following
\begin{proposition}\label{proV3}
For any given $\lambda_3>0$, there exists $\delta_3>0$ and positive constants $K_{10},\ K_{11},\ K_{12}$, $K_{13}$ and $K_{14} $,  such that  
\begin{align}
\dot{V}_3\leq -\lambda_3 V_3+K_{10}\|\gamma\|_{\infty}V_3+K_{11}V_3V_2^{\frac{1}{2}}+K_{12}V_2V_3^{\frac{1}{2}}+K_{13}V_3^{\frac{3}{2}}+K_{14}\|\zeta\|^3_\infty,
\end{align}
provided that $\|\gamma\|_{\infty}+\|\zeta\|_{\infty}\leq \delta_3$.
\end{proposition}
 \subsection{Proof of the  $H^2$ stability  for $\gamma$}
 Denote $W=V_1+V_2+V_3$,  by Proposition \ref{proV1}, \ref{proV2} and \ref{proV3} as well as Lemma \ref{B.7lemma},  one can show that for any given $\lambda>0$, there exists $\delta>0$ and $K_{15}>0$, such that  
 \begin{align}
 \dot{W}\leq -\lambda W+K_{15} W^{\frac{3}{2}},
 \end{align}
 provided that $\|\gamma\|_{\infty}+\|\zeta\|_{\infty}\leq \delta$.
 This concludes the whole proof of Theorem \ref{main resultNL}.\cqfd
 
 \section*{Acknowledgements}
The authors would like to thank Jean-Michel Coron for his encouragement and fruitful discussions, and also thank  Guillaume Olive and Fr\'ed\'eric Marbach for interesting remarks.

\appendix
\renewcommand\thesection{A\Alph{section}}
\renewcommand\theequation{A.\arabic{equation}}
\setcounter{equation}{0}
\setcounter{section}{0}

\section*{Appendix A}\label{appendixA}
In this section, we will show the well-posedness and piecewise smoothness of the Kernel $K$ and $L$ which are given by the following Theorems.
\begin{thm}\label{theorem appendix}
Let $N\in\mathbb{N}^{+}$. Under the assumption that $\sigma_{ij}\in C^{N}[0,1],\ \lambda_i\in C^{N}[0,1](i,j=1,\cdots,n)$, there exists a unique piecewise $C^{N}(\mathcal{T})$ solution $K$ to the  hyperbolic system (\ref{eq:developedKernelK})-(\ref{eq:artificialboundary}). Moreover, if the $C^{N-1}$ compatibility conditions at the point $(x,\xi)=(1,1)$ are satisfied, then $K(\cdot,\cdot)\in C^{N-1}(0,1)$, $K(\cdot,0) \in C^{N-1}(0,1)$ with bounded $C^{N-1}$ norm.

\end{thm}
\textbf{Proof.} We divided the proof into two parts. For the first part, we prove the regularity of the kernels. For this, we only prove the case $N=1$. For $N\geq 1$, the results can be obtained by induction. In the case $N=1$, one can, in fact, refer \cite{hu-dimeglio-vazquez-krstic} and Remark \ref{continuouskernel} to find there exists a piecewise $C^{0}$  kernel $K$
 for the boundary problem (\ref{eq:developedKernelK})-(\ref{eq:artificialboundary}), 
 where though only constant coupling coefficients and transport velocities are considered. However, the method in \cite{hu-dimeglio-vazquez-krstic} straightforwardly extends to spatially varying coefficients with more involved technical developments. Next, we will improve the regulality of $K$.   Let $\mathcal{H}_{ij}=\partial_x K_{ij}(x,\xi)$ and $\mathcal{Y}_{ij}=\partial_\xi K_{ij}(x,\xi)$. By differentiating with respect to $x$  in (\ref{eq:developedKernelK}),  one can show that 
\begin{align}
	\lambda_i(x) \partial_x \mathcal{H}_{ij}(x,\xi)+\lambda_j(\xi) \partial_\xi \mathcal{H}_{ij}(x,\xi) &=- \sum\limits_{k=1}^n\big(\sigma_{kj}(\xi)+\delta_{kj}\lambda_j'(\xi)\big)\mathcal{H}_{ik}(x,\xi)-\lambda'_i(x)\mathcal{H}_{ij}(x,\xi).\label{eq:developedKernelH}
\end{align}
Differentiating the boundary conditions in (\ref{eq:hypotenuseK}) and (\ref{eq:modifx0boundary}), we have
\begin{align}
&\mathcal{H}_{ij}(x,x)+\mathcal{Y}_{ij}(x,x)=k'_{ij}(x) &\text{for } 1\leq i, j\leq n(i\ne j),\label{1}\\
&\mathcal{H}_{ij}(x,0)=-\frac{1}{\lambda_j(0)}\sum\limits_{k=1}^{n-m} \lambda_{m+k}(0) \mathcal{H}_{i,m+k}(x,0)q_{k,j} & \text{for } 1 \leq i \leq j \leq m.\label{boundaryforx=0H}
\end{align}
Next, differentiating the boundary conditions in \eqref{eq:artificialboundary1}--\eqref{eq:artificialboundary}, we have
\begin{align}\label{2}
\mathcal{Y}_{ij}(1,\xi)=\dot{k}^{(1)}_{ij}(\xi), \ \text{for} \ \  1&\leq j < i \leq m \ \cup \ m+1\leq i<j\leq n
\end{align}
and the boundary conditions for $\mathcal{H}_{ij}(n\geq i\geq j\geq m+1)$ on $\xi=0$:
\begin{align}\label{boundaryforx=0Hartificial}
\mathcal{H}_{ij}(x,0)=\dot{k}^{(2)}_{ij}(x), \ \text{for} \ \ m+1&\leq j\leq i \leq n.
\end{align}
In view of the equations (\ref{eq:developedKernelK}), it is easy to see that
\begin{align}
\lambda_i(x) \mathcal{H}_{ij}(x,x)+\lambda_j(x) \mathcal{Y}_{ij}(x,x) &=- \sum\limits_{k=1}^n\big(\sigma_{kj}(x)+\delta_{kj}\lambda_j'(x)\big)K_{ik}(x,x)\label{3}\\
\lambda_i(1) \mathcal{H}_{ij}(1,\xi)+\lambda_j(\xi) \mathcal{Y}_{ij}(1,\xi) &=- \sum\limits_{k=1}^n\big(\sigma_{kj}(\xi)+\delta_{kj}\lambda_j'(\xi)\big)K_{ik}(1,\xi)\label{4}
\end{align}
Combining (\ref{1}) and (\ref{3}), we have
\begin{align}\label{Hxx}
\mathcal{H}_{ij}(x,x)=\frac{\lambda_j(x)k'_{ij}(x)+\sum\limits_{k=1}^n\big(\sigma_{kj}(x)+\delta_{kj}\lambda_j'(x)\big)K_{ik}(x,x)}{\lambda_j(x)-\lambda_i(x)},\ \ \text{for } 1\leq i, j\leq n(i\ne j).
\end{align}
Similarly, plugging (\ref{4}) into (\ref{2}), one immediately obtains, for $1 \leq j < i \leq m \ \cup \ m+1\leq i<j\leq n$, we have
\begin{align}\label{last}
\mathcal{H}_{ij}(1,\xi)=-\frac{1}{\lambda_i(1)} \left(\sum\limits_{k=1}^n\big(\sigma_{kj}(\xi)+\delta_{kj}\lambda_j'(\xi)\big)K_{ik}(1,\xi)+\lambda_j(\xi)\dot{k}^{(1)}_{ij}(\xi)\right),
\end{align}
which are piecewise $C^0(0,1)$ function. By the theory in \cite{hu-dimeglio-vazquez-krstic},  we can prove that there exists a unique piecewise $\mathcal{H}\in C^{0}(\mathcal{T})$ for the boundary value problem (\ref{eq:developedKernelH}),  (\ref{boundaryforx=0H}), (\ref{boundaryforx=0Hartificial}) and (\ref{Hxx})--(\ref{last}). Noting the equations (\ref{eq:developedKernelK}),  we know that $\mathcal{Y}$ shares the same regularity as $\mathcal{H}$.

Next, we prove the regularity of $K(\cdot,0)$. Obviously, for $N=1$, by the theory in \cite{hu-dimeglio-vazquez-krstic} and Remark \ref{continuouskernel}, one can prove that both $K(\cdot,\cdot)$ and $K(\cdot,0)\in C^0(0,1)$ with bounded $C^0$ norm,  provided that the $C^0$ compatibility conditions (\ref{designcomconC0}) are satisfied at the the point $(x,\xi)=(1,1)$.  Next, we prove the case $N=2$. Taking an $\xi$-derivative in (\ref{eq:developedKernelK}) yields
\begin{align}
\begin{aligned}\label{eq:developedKernelY}
	\lambda_i(x) \partial_x \mathcal{Y}_{ij}(x,\xi)+\lambda_j(\xi) \partial_\xi \mathcal{Y}_{ij}(x,\xi) &=- \sum\limits_{k=1}^n\big(\sigma_{kj}(\xi)+\delta_{kj}\lambda_j'(\xi)\big)\mathcal{Y}_{ik}(x,\xi)-\lambda'_j(\xi)\mathcal{Y}_{ij}(x,\xi)\\
	&\ \ \ \ -\sum\limits_{k=1}^n\big(\sigma'_{kj}(\xi)+\delta_{kj}\lambda_j''(\xi)\big)\mathcal{K}_{ik}(x,\xi)
\end{aligned}
\end{align}
Combining (\ref{1}) and (\ref{3}), we have
\begin{align}\label{yxx}
\mathcal{Y}_{ij}(x,x)=\frac{\lambda_i(x)k'_{ij}(x)+\sum\limits_{k=1}^n\big(\sigma_{kj}(x)+\delta_{kj}\lambda_j'(x)\big)K_{ik}(x,x)}{\lambda_i(x)-\lambda_j(x)},\ \ \text{for } 1\leq i, j\leq n(i\ne j).
\end{align}
Since 
\begin{align}\label{x=0yandh}
	\lambda_i(x) \mathcal{H}_{ij}(x,0)+\lambda_j(0) \mathcal{Y}_{ij}(x,0) &=- \sum\limits_{k=1}^n\big(\sigma_{kj}(0)+\delta_{kj}\lambda_j'(0)\big)K_{ik}(x,0)
	\end{align}
Plugging (\ref{boundaryforx=0H}) and (\ref{boundaryforx=0Hartificial}), respectively, one obtains
\begin{align}\label{yx=01}
\mathcal{Y}_{ij}(x,0)=-\frac{1}{\lambda_j(0)}\bigg(\lambda_i(x)\dot{k}^{(2)}_{ij}(x)+\sum\limits_{k=1}^n\big(\sigma_{kj}(0)+\delta_{kj}\lambda_j'(0)\big)K_{ik}(x,0)\bigg), \ \text{for} \ \ m+1&\leq j\leq i \leq n
\end{align}
and 
\begin{align}
\begin{aligned}\label{yx=02}
\mathcal{Y}_{ij}(x,0)=&-\frac{1}{\lambda_j(0)}\sum\limits_{k=1}^n\big(\sigma_{kj}(0)+\delta_{kj}\lambda_j'(0)\big)K_{ik}(x,0)+\frac{1}{\lambda^2_j(0)}\sum_{k=1}^{n-m}\lambda^2_{m+k}(0)q_{k,j}\mathcal{Y}_{i,m+k}(x,0)\\
&+\frac{1}{\lambda^2_j(0)}\sum_{k=1}^{n-m}\sum_{s=1}^n\lambda_{m+k}(0)q_{k,j}\big(\sigma_{s,m+k}(0)+\delta_{s,m+k}\lambda_{m+k}'(0)\big)K_{is}(x,0),\ \text{for } 1 \leq i \leq j \leq m.
\end{aligned}
\end{align}
Noting (\ref{2}), (\ref{yxx}), (\ref{yx=01}) and $K(\cdot,0)\in C^0$, we know that $\mathcal{Y}_{ij}(\cdot,\cdot)\in C^0(0,1)(i\ne j)$. $\mathcal{Y}_{ij}(1,\cdot)\in C^0(0,1)(\text{for}1\leq j < i \leq m \ \cup \ m+1\leq i<j\leq n)$ and $\mathcal{Y}_{ij}(\cdot,0)\in C^0(0,1)(\text{for}\ m+1\leq j\leq i \leq n
 )$. By the $C^1$ compatibility conditions (\ref{designcomconC1}) at the point $(x,\xi)=(1,1)$ and using the theory in \cite{hu-dimeglio-vazquez-krstic} and Remark \ref{continuouskernel}, we can prove that there exists a unique piecewise $C^0$ function $\mathcal{Y}=\mathcal{Y}(x,\xi)$ for the boundary value problem (\ref{eq:developedKernelY}), (\ref{yxx}), (\ref{yx=01}), (\ref{2}) and (\ref{yx=02}),  which satisfies $\mathcal{Y}(\cdot,\cdot),\ \mathcal{Y}(\cdot,0)\in C^{0}(0,1)$. Noting (\ref{x=0yandh}) and (\ref{3}), we know that $\mathcal{H}(\cdot,\cdot),\ \mathcal{H}(\cdot,0)\in C^{0}(0,1)$. This finishes the proof.\cqfd

\begin{rem}\label{continuouskernel}
It is worthy of mentioning that in \cite{hu-dimeglio-vazquez-krstic}, we only prove $K\in L^{\infty}(\mathcal{T})$ and do not clarify the regularity of the kernel because of  brevity purposes. However, with the same procedure in \cite[Section A.3]{JMC-Vazquez} and \cite{Meglio}, one can prove that $K$ is a piecewise $C^0$ function with $K(\cdot,\cdot),\ K(\cdot,0)\in C^{0}(0,1)$ and $K(1,\cdot)$ being a function of {piecewise} $C^0(0,1)$ for the boundary problem (\ref{eq:developedKernelK})-(\ref{eq:artificialboundary}), provided $\sigma_{ij}\in C^{0}[0,1],\ \lambda_i\in C^{1}[0,1](i,j=1,\cdots,n)$ and the $C^0$ compatibility conditions (\ref{designcomconC0}) are satisfied at the the point $(x,\xi)=(1,1)$.
\end{rem} 

\begin{thm}\label{theorem appendix2}
Under the assumptions of Theorem \ref{theorem appendix}, For any $N\in\mathbb{N}$, there exists a unique piecewise $C^{N}(\mathcal{T})$ kernel $L$ to the inverse transformation (\ref{4.1}). Moreover, $L(x,x),\ L(x,0)\in C^{N-1}(0,1)$.
\end{thm}
\textbf{Proof.} 
Substituting (\ref{2.1}) for (\ref{4.1}), it is easy to see that $L$ is the solution of the following Volterra equations
\begin{align}\label{inverseL}
L(x,\xi)=K(x,\xi)+\int_\xi^xK(x,s)L(s,\xi)ds,
\end{align}
which yields that
\begin{align}
L(x,x)=K(x,x) \in C^{N-1}(0,1).
\end{align}
Noting (\ref{boundaryKxx}), we have
\begin{align}\label{invLxx}
\Sigma(x)+L(x,x)\Lambda(x)-\Lambda(x)L(x,x)=0.
\end{align}
Next, Taking a partial derivative in $x$ and $\xi$ in (\ref{inverseL}), respectively, one obtains
\begin{align}
&L_x(x,\xi)=K_x(x,\xi)+K(x,x)L(x,\xi)+\int_\xi^xK_x(x,\xi)L(s,\xi)ds,\label{L_x}\\
&L_\xi(x,\xi)=K_\xi(x,\xi)-K(x,\xi)L(\xi,\xi)+\int_\xi^xK(x,s)L_\xi(s,\xi)ds.\label{L_xi}
\end{align}
Substituting (\ref{L_x}) and (\ref{L_xi}) for (\ref{2.6}) and using integration by parts, one has
\begin{align}\label{invkernelequation}
\Lambda(x)L_x(x,\xi)+L_{\xi}(x,\xi)\Lambda(\xi)=\big(\Sigma(x)-\Lambda_\xi(\xi)\big)L(x,\xi)
\end{align}
Again by (\ref{inverseL}), we have
\begin{align}\label{Lx=0}
L(x,0)=K(x,0)+\int_0^x K(x,s)L(s,0)ds.
\end{align}
Since both $K(x,0)$ and $K(x,x)$ are $C^{N-1}$ continuous functions, by a suitable iteration procedure (see \cite[Theorem 3.2, Pages 32--34]{Linz-1985}), it easy to see that there exists $L(x,0)=l(x)\in C^{N-1}(0,1)$ for the Volterra equation of the second kind (\ref{Lx=0}).
 
On the other hand, substituting (\ref{4.1}) for (\ref{2.1}), one gets
\begin{align}
L(x,\xi)=K(x,\xi)+\int_\xi^x L(x,s)K(s,\xi)ds,
\end{align}
then
\begin{align}
L(1,\xi)=K(1,\xi)+\int_\xi^1 L(1,s)K(s,\xi)ds.
\end{align}
With the same argument above, we can see that $L_{ij}(1,\xi)=\tilde{l}_{ij}(\xi)(m\geq i>j\geq 1, n\geq j>i\geq m+1)$ on $x=1$ are functions of piecewise $C^{N}(0,1)$. Then, For the boundary problem (\ref{invkernelequation}) with the boundary conditions (\ref{invLxx}), and
\begin{align}
L_{ij}(1,\xi)=\tilde{l}_{ij}(\xi), \ \text{for} \ \ 1&\leq j < i \leq m \ \cup \ m+1\leq i<j\leq n\\
L_{ij}(x,0)=l_{ij}(x),\ \text{for}\ \ 1&\leq i\leq j\leq m\ \cup\ m+1\leq j\leq i\leq n.
\end{align}
by  Theorem \ref{theorem appendix}, one immediately gets Theorem \ref{theorem appendix2}.\cqfd

\renewcommand\thesection{B\Alph{section}}
\renewcommand\theequation{B.\arabic{equation}}
\setcounter{equation}{0}
\setcounter{lemma}{0}
\section*{Appendix B}\label{appendixB}
In this appendix, we first sketch out four useful lemmas (the details can be found in \cite{JMC-Vazquez}).   

\begin{lemma}\label{b.1}
There exists a positive real number $c_1$, such that
\begin{align}
|\mathcal{K}[\gamma]|+|\mathcal{L}[\gamma]|+|\mathcal{K}_1[\gamma]|+|\mathcal{K}_2[\gamma]|+|\mathcal{L}_1[\gamma]|\leq c_1(|\gamma|+\|\gamma\|_{L^1}).
\end{align}
\end{lemma}
\begin{lemma}\label{B.2lemma}
Suppose $\|\gamma\|_{\infty}$ is suitable small, one can see that
\begin{align}
&|F_1|\leq c_2(|\gamma|+\|\gamma\|_{L^1}),\\
&|F_2|\leq c_3(|\gamma|^2+\|\gamma\|^2_{L^1}),\\
&|F_3|\leq c_4(|\gamma|+\|\gamma\|_{L^1})(\|\gamma_x\|_{L^2}+|\gamma_x|),\\
&|F_4|\leq c_5(|\gamma|^2+\|\gamma\|^2_{L^1}).
\end{align}
\end{lemma}

\begin{lemma}\label{B.3lemma}
\begin{align}
&|F_{11}|\leq c_6(|\zeta|+\|\zeta\|_{L^1}),\\
&|F_{12}|\leq c_7(|\gamma_x|+|\gamma|+\|\gamma\|_{L^1}),\\
&|F_{21}|\leq c_8(|\gamma|+\|\gamma\|_{L^1})(|\zeta|+\|\zeta\|_{L^1}),\\
&|F_{5}|\leq c_9(|\zeta|+\|\zeta\|_{L^2})(|\gamma|+\|\gamma\|_{L^2})+c_{10}(|\zeta|+\|\zeta\|_{L^2})((|\gamma_x|+\|\gamma_x\|_{L^2}))+c_{11}|\gamma(0)||\zeta(0)|,\\
&|F_6|\leq c_{12}(|\gamma|+\|\gamma\|_{L^2})(|\zeta|+\|\zeta\|_{L^2}).
\end{align}
\end{lemma}

\begin{lemma}
\begin{align}
|F_{13}|\leq& c_{13}(|\zeta|^2+\|\zeta\|^2_{L^2})+c_{14}(|\theta|+\|\theta\|_{L^1}),\\
|F_{14}|\leq& c_{14}(|\zeta|+\|\zeta\|_{L^1})(1+|\gamma_x|+|\gamma|+\|\gamma\|_{L^1})+c_{15}(|\zeta|+|\zeta_x|+\|\zeta\|_{L^1}),\\
|F_{22}|\leq& c_{16}(|\gamma|+\|\gamma\|_{L^1})(|\theta|+\|\theta\|_{L^1})+c_{17}(|\zeta|^2+\|\zeta\|^2_{L^2}),\\
|F_7|\leq& c_{18}(|\zeta|^2+\|\zeta\|^2_{L^2})(1+|\gamma|+\|\gamma_x\|)\notag\\
 &+c_{19}(|\zeta|+\|\zeta\|_{L^2})(|\zeta_x|+\|\zeta\|_{L^2})\\
 &+c_{20}(|\gamma|+\|\gamma\|_{L^2}+|\gamma_x|)(|\theta|+\|\theta\|_{L^2})\notag\\
 &+c_{21}(|\zeta(0)|^2+|\gamma(0)||\theta(0)|),\notag\\
 |F_8|\leq& c_{22}(|\zeta|^2+\|\zeta\|^2_{L^2})(1+\|\gamma\|_{{\infty}})+c_{23}(|\gamma|+\|\gamma\|_{L^2})(|\theta|+\|\theta\|_{L^2}).
\end{align}
\end{lemma}

Next, we show the following proposition which is also mentioned in \cite{JMC-Vazquez}, however here  more technical developments are involved.
\setcounter{proposition}{0}
\begin{proposition}\label{appendixb-pro}
There exists $\delta>0$ such that for any $|\gamma|+|\zeta|\leq \delta$, one has
\begin{align}
\|\theta\|_{\infty}&\leq C_1(\|\gamma_{xx}\|_{\infty}+\|\gamma_x\|_{\infty}+\|\gamma\|_{\infty}),\\
\|\theta\|_{L^2}&\leq C_2(\|\gamma_{xx}\|_{L^2}+\|\gamma_x\|_{L^2}+\|\gamma\|_{L^2}),\\
\|\gamma_{xx}\|_{\infty}&\leq C_3(\|\theta\|_{\infty}+\|\zeta\|_{\infty}+\|\gamma\|_{\infty}),\\
\|\gamma_{xx}\|_{L^2}&\leq C_4(\|\theta\|_{L^2}+\|\zeta\|_{L^2}+\|\gamma\|_{L^2}),
\end{align}
where $C_1$, $C_2$, $C_3$ and $C_4$ are positive constants.
\end{proposition}
\textbf{Proof.} We prove the next  three lemmas to get Proposition \ref{appendixb-pro}.
\begin{lemma}\label{B.5lemma}
There exists $\delta$ such that, if $|\gamma|\leq \delta$, then the following inequalities hold:
\begin{align}
\|\zeta\|_{{\infty}}&\leq c_1(\|\gamma_x\|_{\infty}+\|\gamma\|_{\infty})\label{b.20}\\
\|\zeta\|_{L^{2}}&\leq c_2(\|\gamma_x\|_{L^2}+\|\gamma\|_{L^2}),\label{b.21}\\
\|\gamma_x\|_{{\infty}}&\leq c_3(\|\zeta\|_{\infty}+\|\gamma\|_{\infty}),\label{b.22}\\
\|\gamma_x\|_{L^{2}}&\leq c_4(\|\zeta\|_{L^2}+\|\gamma\|_{L^2})\label{b.23}
\end{align}
\end{lemma}
\textbf{Proof.} Noting (\ref{analysis gamma}), one can easily see that
\begin{align}\label{b.24}
\begin{aligned}
\zeta(t,x)+\Lambda(x)\gamma_x(t,x)-G(x)\gamma(t,0)=F_3[\gamma,\gamma_x]+F_4[\gamma].
\end{aligned}
\end{align}
The difference between our proof and the proof in \cite[Lemma B.6]{JMC-Vazquez} is the appearance of the term $G(x)\gamma(t,0)$ in (\ref{b.24}). Noting (\ref{2.29}) and Theorem \ref{theorem appendix}, we have $G(\cdot)\in C^1(0,1)$ with bounded $C^1$ norm. Then since one can show that
\begin{align}
\|G(\cdot)\gamma(t,0)\|_{L^2}\leq c_5\|G(\cdot)\gamma(t,0)\|_{\infty}\leq c_6\|\gamma\|_{\infty}\leq c_7(\|\gamma_x\|_{L^2}+\|\gamma\|_{L^2}),
\end{align}
which yields, by the same argument in \cite[Lemma B.6]{JMC-Vazquez}, (\ref{b.20})-(\ref{b.22}).

On the other hand, by the special structure of $G(x)$, we have
\begin{align}
&\|\partial_x\gamma_1\|_{L^2}\leq c_8(\|\zeta\|_{L^2}+\|\gamma_x\|_{L^2}\|\gamma\|_{{\infty}}+\|\gamma\|_{L^2}\|\gamma\|_{{\infty}}),\\
&\|\partial_x\gamma_2\|_{L^2}\leq c_9(\|\zeta\|_{L^2}+\|\gamma_1\|_{\infty}+\|\gamma_x\|_{L^2}\|\gamma\|_{{\infty}}+\|\gamma\|_{L^2}\|\gamma\|_{{\infty}}),\\
&\vdots\ \ \ \ \ \ \ \ \ \ \ \ \ \ \ \ \ \ \ \ \ \  \ \vdots\ \ \ \ \ \ \ \ \ \ \ \ \ \ \ \ \ \ \ \ \ \  \ \vdots\ \ \ \ \ \ \ \ \ \ \ \ \ \ \ \ \ \ \ \ \ \  \ \vdots\notag\\
&\|\partial_x\gamma_m\|_{L^2}\leq c_{m+7}(\|\zeta\|_{L^2}+\sum_{r=1}^{m-1}\|\gamma_r\|_{\infty}+\|\gamma_x\|_{L^2}\|\gamma\|_{{\infty}}+\|\gamma\|_{L^2}\|\gamma\|_{{\infty}}),\\
&\|\partial_x\gamma_s\|_{L^2}\leq c_{s+7}(\|\zeta\|_{L^2}+\sum_{r=1}^{m}\|\gamma_r\|_{\infty}+\|\gamma_x\|_{L^2}\|\gamma\|_{{\infty}}+\|\gamma\|_{L^2}\|\gamma\|_{{\infty}}),
\end{align}
in which $s=m+1,\cdots,n$. Noting the classical Sobolev's inequality
 \begin{align}
   &\|\gamma\|_{L^{\infty}}\leq \widetilde{C}\Big(\|\gamma\|_{L^2}+\|\gamma_x\|_{L^2}\Big)\leq \widetilde{\overline{C}}\|\gamma\|_{H^1},\label{sobolev2}
 \end{align}
one gets that
\begin{align}
&\|\partial_x\gamma_1\|_{L^2}\leq C_1(\|\zeta\|_{L^2}+\|\gamma_x\|_{L^2}\|\gamma\|_{{\infty}}+\|\gamma\|_{L^2}\|\gamma\|_{{\infty}}),\\
&\|\partial_x\gamma_2\|_{L^2}\leq C_2(\|\zeta\|_{L^2}+\|\gamma\|_{L^2}+\|\partial_x\gamma_1\|_{L^2}+\|\gamma_x\|_{L^2}\|\gamma\|_{{\infty}}+\|\gamma\|_{L^2}\|\gamma\|_{{\infty}}),\\
&\vdots\ \ \ \ \ \ \ \ \ \ \ \ \ \ \ \ \ \ \ \ \ \  \ \vdots\ \ \ \ \ \ \ \ \ \ \ \ \ \ \ \ \ \ \ \ \ \  \ \vdots\ \ \ \ \ \ \ \ \ \ \ \ \ \ \ \ \ \ \ \ \ \  \ \vdots\notag\\
&\|\partial_x\gamma_m\|_{L^2}\leq C_{m}(\|\zeta\|_{L^2}+\|\gamma\|_{L^2}+\sum_{r=1}^{m-1}\|\gamma_r\|_{L^2}+\|\gamma_x\|_{L^2}\|\gamma\|_{{\infty}}+\|\gamma\|_{L^2}\|\gamma\|_{{\infty}}),\\
&\|\partial_x\gamma_s\|_{L^2}\leq C_{s}(\|\zeta\|_{L^2}+\|\gamma\|_{L^2}+\sum_{r=1}^{m}\|\gamma_r\|_{L^2}+\|\gamma_x\|_{L^2}\|\gamma\|_{L^{\infty}}+\|\gamma\|_{L^2}\|\gamma\|_{{\infty}}),
\end{align}
where $s=m+1,\cdots,n$.
Then, we can easily obtain by induction that 
\begin{align}
\|\gamma_x\|_{L^2}\leq \tilde{c}_1(\|\zeta\|_{L^2}+\|\gamma_x\|_{L^2}\|\gamma\|_{{\infty}}+\|\gamma\|_{L^2}\|\gamma\|_{{\infty}}+\|\gamma\|_{L^2}),
\end{align}
which concludes (\ref{b.23}), under the assumption that $\|\gamma\|_{\infty}$ is small enough. \cqfd

Combining the same technical approach as in \cite[Lemma B.7 and Lemma B.8]{JMC-Vazquez} and an analogous argument used in the proof of Lemma \ref{B.5lemma} and noting $G\in C^1$, the details of which we omit, one can show the next two lemmas.
\begin{lemma}\label{B.6lemma}
There exists $\delta$ such that, if $\|\gamma\|_{{\infty}}\leq \delta$, then the following inequalities hold:
\begin{align}
&\|\gamma_{xx}\|_{\infty}\leq c_1(\|\zeta_x\|_{{\infty}}+\|\zeta\|_{{\infty}}+\|\gamma\|_{\infty}),\\
&\|\gamma_{xx}\|_{L^2}\leq c_2(\|\zeta_x\|_{L^{2}}+\|\zeta\|_{L^{2}}+\|\gamma\|_{L^2}),\\
&\|\zeta_x\|_{\infty}\leq c_3(\|\gamma_{xx}\|_{{\infty}}+\|\zeta\|_{{\infty}}+\|\gamma\|_{\infty}),\\
&\|\zeta_x\|_{L^2}\leq c_4(\|\gamma_{xx}\|_{L^{2}}+\|\zeta\|_{L^{2}}+\|\gamma\|_{L^2}),
\end{align}
where $c_1$, $c_2$, $c_3$ and $c_4$ are positive constants.
\end{lemma}
and
\begin{lemma}\label{B.7lemma}
There exists $\delta$ such that, if $\|\gamma\|_{{\infty}}+\|\zeta\|_{\infty}\leq \delta$, then the following inequalities hold:
\begin{align}
&\|\theta\|_{\infty}\leq c_1(\|\zeta_x\|_{{\infty}}+\|\zeta\|_{{\infty}}+\|\gamma\|_{\infty}),\\
&\|\theta\|_{L^2}\leq c_2(\|\zeta_x\|_{L^{2}}+\|\zeta\|_{L^{2}}+\|\gamma\|_{L^2}),\\
&\|\zeta_x\|_{\infty}\leq c_3(\|\theta\|_{{\infty}}+\|\zeta\|_{{\infty}}+\|\gamma\|_{\infty}),\\
&\|\zeta_x\|_{L^2}\leq c_4(\|\theta\|_{L^{2}}+\|\zeta\|_{L^{2}}+\|\gamma\|_{L^2}),
\end{align}
where $c_1$, $c_2$, $c_3$ and $c_4$ are positive constants.
\end{lemma}
The above three Lemma \ref{B.5lemma}--\ref{B.7lemma} immediately yield  Proposition \ref{appendixb-pro}. \cqfd

\bibliographystyle{plain}
\bibliography{biblio}

\end{document}